# Doubly Accelerated Methods for
# Faster CCA and Generalized Eigendecomposition

Zeyuan Allen-Zhu
zeyuan@csail.mit.edu
Institute for Advanced Study

Yuanzhi Li
yuanzhil@cs.princeton.edu
Princeton University

July 20, 2016*

### Abstract

We study $k$-GenEV, the problem of finding the top $k$ generalized eigenvectors, and $k$-CCA, the problem of finding the top $k$ vectors in canonical-correlation analysis. We propose algorithms `LazyEV` and `LazyCCA` to solve the two problems with running times linearly dependent on the input size and on $k$.

Furthermore, our algorithms are *doubly-accelerated*: our running times depend only on the square root of the matrix condition number, and on the square root of the eigengap. This is the first such result for both $k$-GenEV or $k$-CCA. We also provide the first gap-free results, which provide running times that depend on $1/\sqrt{\varepsilon}$ rather than the eigengap.

## 1 Introduction

The Generalized Eigenvector (GenEV) problem and the Canonical Correlation Analysis (CCA) are two fundamental problems in scientific computing, machine learning, operations research, and statistics. Algorithms solving these problems are often used to extract features to compare large-scale datasets, as well as used for problems in regression [17], clustering [10], classification [18], word embeddings [11], and many others.

Given two symmetric matrices $A, B \in \mathbb{R}^{d \times d}$ where $B$ is positive definite. The GenEV problem is to find *generalized eigenvectors* $v_1, \dots, v_d$ where each $v_i$ satisfies

$$v_i \in \arg\max_{v \in \mathbb{R}^d} |v^\top A v| \quad \text{such that} \quad \left\{ \begin{array}{l} v^\top B v = 1 \\ v^\top B v_j = 0 \ \forall j \in [i-1] \end{array} \right.$$

The values $\lambda_i \stackrel{\text{def}}{=} v_i^\top A v_i$ are known as the generalized eigenvalues, and it satisfies $|\lambda_1| \geq \cdots |\lambda_d|$. Following the tradition of [13, 29], we assume *without loss of generality* that $\lambda_i \in [-1, 1]$.

Given matrices $X \in \mathbb{R}^{n \times d_x}, Y \in \mathbb{R}^{n \times d_y}$ and denoting by $S_{xx} = \frac{1}{n} X^\top X$, $S_{xy} = \frac{1}{n} X^\top Y$, $S_{yy} = \frac{1}{n} Y^\top Y$, the CCA problem is to find *canonical-correlation vectors* $\{(\phi_i, \psi_i)\}_{i=1}^r$ where $r = \min\{d_x, d_y\}$ and each pair

$$(\phi_i, \psi_i) \in \arg\max_{\phi \in \mathbb{R}^{d_x}, \psi \in \mathbb{R}^{d_y}} \left\{ \phi^\top S_{xy} \psi \right\} \quad \text{such that} \quad \left\{ \begin{array}{l} \phi^\top S_{xx} \phi = 1 \wedge \phi^\top S_{xx} \phi_j = 0 \ \forall j \in [i-1] \\ \psi^\top S_{yy} \psi = 1 \wedge \psi^\top S_{yy} \psi_j = 0 \ \forall j \in [i-1] \end{array} \right.$$

The values $\sigma_i \stackrel{\text{def}}{=} \phi_i^\top S_{xy} \psi_i \geq 0$ are known as the canonical-correlation coefficients, and it *always* satisfies $1 \geq \sigma_1 \geq \cdots \geq \sigma_r \geq 0$.

---

*First appeared on arXiv on this date. In this new version, we have stated more clearly why this paper has outperformed relevant previous results, and included discussions for doubly-stochastic methods.

It is a folklore that solving CCA *exactly* can be reduced to solving GenEV exactly, if one defines $B = \text{diag}\{S_{xx}, S_{yy}\} \in \mathbb{R}^{d \times d}$ and $A = [[0, S_{xy}]; [S_{xy}^\top, 0]] \in \mathbb{R}^{d \times d}$ for $d \stackrel{\text{def}}{=} d_x + d_y$, see Lemma 2.3.

Despite the fundamental importance and the frequent necessity in applications, there are few results on obtaining provably efficient algorithms for GenEV and CCA until very recently. In the breakthrough result of Ma, Lu and Foster [20], they proposed to study algorithms to find top $k$ generalized eigenvectors ($k$-GenEV) or top $k$ canonical-correlation vectors ($k$-CCA). They designed an alternating minimization algorithm whose running time is only linear in terms of the number of non-zero elements of the matrix (that we denote by $\text{nnz}(A)$ for a matrix $A$ in this paper), and also nearly-linear in $k$. Such algorithms are very appealing because in real-life applications, it is often only relevant to obtain top correlation vectors, as opposed to the less meaningful vectors in the directions where the datasets do not correlate. Unfortunately, the method of Ma, Lu and Foster has a running time that linearly scales with $\kappa$ and $1/\text{gap}$, where

- $\kappa \geq 1$ is the condition number of matrix $B$ in GenEV, or of matrices $X^\top X, Y^\top Y$ in CCA; and
- $\text{gap} \in [0, 1)$ is the (relative) eigengap $\frac{\lambda_k - \lambda_{k+1}}{\lambda_k}$ in GenEV, or $\frac{\sigma_k - \sigma_{k+1}}{\sigma_k}$ in CCA.

These parameters are usually not constants and scale with the problem size.

### Challenge 1: Acceleration

For many easier scientific computing problems, we are able to design algorithms that have *accelerated* dependencies on $\kappa$ and $1/\text{gap}$. As two concrete examples, $k$-PCA can be solved with a running time linearly in $1/\sqrt{\text{gap}}$ as opposed to $1/\text{gap}$ [16]; computing $B^{-1}w$ for a vector $w$ can be solved in time linearly in $\sqrt{\kappa}$ as opposed to $\kappa$, where $\kappa$ is the condition number of matrix $B$ [9, 23, 27].

Therefore, can we obtain ***doubly-accelerated*** methods for $k$-GenEV and $k$-CCA, meaning that the running times linearly scale with both $\sqrt{\kappa}$ and $1/\sqrt{\text{gap}}$? Before this paper, for the general case $k > 1$, the method of Ge *et al.* [15] made acceleration possible for parameter $\kappa$, but not for parameter $1/\text{gap}$ (see Table 2).

### Challenge 2: Gap-Freeness

Since $\text{gap}$ can be even zero in the extreme case, can we design algorithms that do not scale with $1/\text{gap}$? Recall that this is possible for the easier task of $k$-PCA. The block Krylov method [22] runs in time linear in $1/\sqrt{\varepsilon}$ as opposed to $1/\sqrt{\text{gap}}$, where $\varepsilon$ is the approximation ratio.

There is no gap-free result previously known for $k$-GenEV or $k$-CCA even for $k = 1$.

### Challenge 3: Stochasticity

For matrix-related problems, one can often obtain a stochastic running time which requires some notations to describe. Consider the simple task of computing $B^{-1}w$ for some vector $w$, where recall accelerated methods solve this problem with a running time linearly in $\sqrt{\kappa}$ for $\kappa$ being the condition number of $B$. If $B = \frac{1}{n} X^\top X$ is given as the covariance matrix of $X \in \mathbb{R}^{n \times d}$, then (accelerated) stochastic gradient methods can be used to compute $B^{-1}w$ in a time linearly in $(1 + \sqrt{\kappa'/n})$ instead of $\sqrt{\kappa}$, where $\kappa' = \frac{\max_{i \in [n]}\{\|X_i\|^2\}}{\lambda_{\min}(B)} \in [\kappa, n\kappa]$ and $X_i$ is the $i$-th row of $X$. (See Lemma 2.6.) Since $1 + \sqrt{\kappa'/n} \leq O(\sqrt{\kappa})$, stochastic methods are *no slower than* non-stochastic ones.

Therefore, can we obtain a similar (but doubly-accelerated!) stochastic method for $k$-CCA?[1] Note that, if the doubly-accelerated requirement is dropped, this task is easier and indeed possible, see Ge *et al.* [15]. However, since their stochastic method is not doubly-accelerated, in certain parameter regimes, it can run slower than non-stochastic ones (even for $k = 1$, see Table 1).

---

[1] Note that a similar problem can be also asked for $k$-GenEV when $A$ and $B$ are both given in their covariance matrix forms. We refrain from doing it in this paper for notational simplicity.

**Remark.** In general, if designed properly,

- Accelerated results are better because they are never slower than non-accelerated ones.
- Gap-free results are better because they imply gap-dependent ones.[2]
- Stochastic results are better because they are never slower than non-stochastic ones.

## 1.1 Our Main Results

We provide algorithms `LazyEV` and `LazyCCA` that are *doubly-accelerated*, *gap-free*, and *stochastic*.[3]

For the general $k$-GenEV problem, our `LazyEV` can be implemented to run in time

$$\widetilde{O}\Big(\frac{k\,\mathsf{nnz}(B)\sqrt{\kappa}}{\sqrt{\mathsf{gap}}} + \frac{k\,\mathsf{nnz}(A) + k^2 d}{\sqrt{\mathsf{gap}}}\Big) \quad \text{or} \quad \widetilde{O}\Big(\frac{k\,\mathsf{nnz}(B)\sqrt{\kappa}}{\sqrt{\varepsilon}} + \frac{k\,\mathsf{nnz}(A) + k^2 d}{\sqrt{\varepsilon}}\Big)$$

in the gap-dependent and gap-free cases respectively. Since our running time only linearly depends on $\sqrt{\kappa}$ and $\sqrt{\mathsf{gap}}$ (resp. $\sqrt{\varepsilon}$), our algorithm `LazyEV` is *doubly-accelerated*.

For the general $k$-CCA problem, our `LazyCCA` can be implemented to run in time

$$\widetilde{O}\Big(\frac{k\,\mathsf{nnz}(X,Y)\cdot\big(1+\sqrt{\kappa'/n}\big) + k^2 d}{\sqrt{\mathsf{gap}}}\Big) \quad \text{or} \quad \widetilde{O}\Big(\frac{k\,\mathsf{nnz}(X,Y)\cdot\big(1+\sqrt{\kappa'/n}\big) + k^2 d}{\sqrt{\varepsilon}}\Big)$$

in the gap-dependent and gap-free cases respectively. Here, $\mathsf{nnz}(X,Y) = \mathsf{nnz}(X) + \mathsf{nnz}(Y)$ and $\kappa' = \frac{2\max_i\{\|X_i\|^2, \|Y_i\|^2\}}{\lambda_{\min}(\mathrm{diag}\{S_{xx}, S_{yy}\})}$ where $X_i$ or $Y_i$ is the $i$-th row vector of $X$ or $Y$. Therefore, our algorithm `LazyCCA` is *doubly-accelerated* and *stochastic*.

We fully compare our running time with prior work in Table 1 and Table 2, and summarize our main contributions below.

- For $k > 1$, we have outperformed all relevant prior works (see Table 2). Moreover, no known method was doubly-accelerated even in the non-stochastic setting.
- For $k \geq 1$, we have obtained the first gap-free running times.
- Even for $k = 1$, we have outperformed most of the state-of-the-arts (see Table 1).

**Other Contributions.** Besides the aforementioned running time improvements, we summarize some other virtues of our algorithms as follows:

- For GenEV, our `LazyEV` distinguishes positive generalized eigenvalues from negative ones. For instance, if $A$ has two generalized eigenvectors $v_1, v_2$ with respect to $B$, one with eigenvalue $\lambda$ and the other with $-\lambda$. Then, previous result GenELin only finds the subspace spanned by $v_1, v_2$ but cannot distinguish $v_1$ from $v_2$.
- For CCA with $k > 1$, previous result CCALin only outputs the subspace spanned by the top $k$ correlation vectors but not identify which vector gives the highest correlation and so on. Our `LazyCCA` provides per-vector guarantees on all the top $k$ correlation vectors, see Theorem 6.3.
- Our `LazyEV` and `LazyCCA` reduce the underlying non-convex problem to multiple calls of quadratic minimization. Since quadratic minimization is a well-studied convex optimization problem, many efficient and robust algorithms can be found. In contrast, previous results for the $k > 1$ case rely on more sophisticated nonconvex optimization; and the previous work of [29] —although uses convex optimization to solve 1-CCA— requires one to work with a sum-of-non-convex function which is usually less efficient to minimize.

---

[2] If a method depends on $1/\varepsilon$ then one can choose $\varepsilon = \mathsf{gap}$ and this translates to a gap-dependent running time.

[3] Recalling Footnote 1, for notational simplicity, we only state our $k$-GenEV result in non-stochastic running time.

Table 1:

| **Problem** | **Paper** | **Running time**   (× for beaten) | | **gap-free?** | **negative EV?** |
|---|---|---|---|---|---|
| 1-GenEV | GenELin [15] | $\widetilde{O}\big(\frac{\mathsf{nnz}(B)\sqrt{\kappa_B}}{\mathsf{gap}} + \frac{\mathsf{nnz}(A)}{\mathsf{gap}}\big)$ | × | no | no |
| | LazyEV **Theorem 4.3** | $\widetilde{O}\big(\frac{\mathsf{nnz}(B)\sqrt{\kappa_B}}{\sqrt{\mathsf{gap}}} + \frac{\mathsf{nnz}(A)}{\sqrt{\mathsf{gap}}}\big)$ | | no | yes |
| | LazyEV **Theorem 4.4** | $\widetilde{O}\big(\frac{\mathsf{nnz}(B)\sqrt{\kappa_B}}{\sqrt{\varepsilon}} + \frac{\mathsf{nnz}(A)}{\sqrt{\varepsilon}}\big)$ | | yes | yes |

| **Problem** | **Paper** | **Running time**   (× for beaten) | | **gap-free?** | **stochastic?** |
|---|---|---|---|---|---|
| 1-CCA | AppGrad [20] | $\mathsf{nnz}(X,Y) \cdot \widetilde{O}\big(\frac{\kappa}{\mathsf{gap}}\big)$ | × | no | no |
| | CCALin [15] | $\mathsf{nnz}(X,Y) \cdot \widetilde{O}\big(\frac{\sqrt{\kappa}}{\mathsf{gap}}\big)$ | × | no | no |
| | ALS [29] | $\mathsf{nnz}(X,Y) \cdot \widetilde{O}\big(\frac{\sqrt{\kappa}}{\mathsf{gap}^2}\big)$ | × | no | no |
| | SI [29] | $\mathsf{nnz}(X,Y) \cdot \widetilde{O}\big(\frac{\sqrt{\kappa}}{\sqrt{\mathsf{gap}}\cdot\sigma_1}\big)$ | × | no | no |
| | CCALin [15] | $\mathsf{nnz}(X,Y) \cdot \widetilde{O}\big(\frac{1+\sqrt{\kappa'/n}}{\mathsf{gap}}\big)$ | × | no | yes |
| | ALS [29] | $\mathsf{nnz}(X,Y) \cdot \widetilde{O}\big(\frac{1+\sqrt{\kappa'/n}}{\mathsf{gap}^2}\big)$ | × | no | yes |
| | LazyCCA **Theorem 6.2** | $\mathsf{nnz}(X,Y) \cdot \widetilde{O}\big(\frac{1+\sqrt{\kappa'/n}}{\sqrt{\mathsf{gap}}}\big)$ | | no | yes |
| | LazyCCA **Theorem 6.3** | $\mathsf{nnz}(X,Y) \cdot \widetilde{O}\big(\frac{1+\sqrt{\kappa'/n}}{\sqrt{\varepsilon}}\big)$ | | yes | yes |
| | SI [29] | $\mathsf{nnz}(X,Y) \cdot \widetilde{O}\big(1 + \frac{\sqrt{\kappa'}/n^{1/4}}{\sqrt{\mathsf{gap}}\cdot\sigma_1}\big)$ *(see Remark 3)* | | no | doubly |
| | LazyCCA **Theorem 6.2** | $\mathsf{nnz}(X,Y) \cdot \widetilde{O}\big(1 + \frac{\sqrt{\kappa'}/n^{1/4}}{\sqrt{\mathsf{gap}}\cdot\sigma_1}\big)$ | | no | doubly |
| | LazyCCA **Theorem 6.3** | $\mathsf{nnz}(X,Y) \cdot \widetilde{O}\big(1 + \frac{\sqrt{\kappa'}/n^{1/4}}{\sqrt{\varepsilon}\cdot\sigma_1}\big)$ | | yes | doubly |

Table 1: Performance comparison on 1-GenEV and 1-CCA.
In GenEV, $\mathsf{gap} = \frac{\lambda_1-\lambda_2}{\lambda_1} \in [0,1]$ and $\kappa_B = \frac{\lambda_{\max}(B)}{\lambda_{\min}(B)} > 1$.
In CCA, $\mathsf{gap} = \frac{\sigma_1-\sigma_2}{\sigma_1} \in [0,1]$, $\kappa = \frac{\lambda_{\max}(\mathrm{diag}\{S_{xx},S_{yy}\})}{\lambda_{\min}(\mathrm{diag}\{S_{xx},S_{yy}\})} > 1$, $\kappa' = \frac{2\max_i\{\|X_i\|^2, \|Y_i\|^2\}}{\lambda_{\min}(\mathrm{diag}\{S_{xx},S_{yy}\})} \in [\kappa, 2n\kappa]$, and $\sigma_1 \in [0,1]$.

**Remark 1.** Stochastic methods depend on modified condition number $\kappa'$; the reason $\kappa' \in [\kappa, 2n\kappa]$ is in Def. 2.4.

**Remark 2.** All non-stochastic CCA methods in this table have been outperformed because $1 + \sqrt{\kappa'/n} \le O(\kappa)$.

**Remark 3.** Doubly-stochastic methods are not necessarily interesting. We discuss them in Section 1.2.

**Remark 4.** Some CCA methods have a running time dependency on $\sigma_1 \in [0,1]$, and this is intrinsic and *cannot* be removed. In particular, if we scale the data matrix $X$ and $Y$, the value $\sigma_1$ stays the same.

**Remark 5.** The only (non-doubly-stochastic) doubly-accelerated method before our work is SI [29] (for 1-CCA only). Our LazyEV is faster than theirs by a factor $\Omega(\sqrt{n\kappa/\kappa'} \times \sqrt{1/\sigma_1})$. Here, $n\kappa/\kappa' \ge 1/2$ and $1/\sigma_1 \ge 1$ are two scaling-invariant quantities usually much greater than 1.

| Problem | Paper | Running time     (✗ for beaten) | gap-free? | negative EV? |
|---|---|---|---|---|
| | GenELin [15] | $\widetilde{O}\left(\frac{knnz(B)\sqrt{\kappa_B}}{gap} + \frac{knnz(A)+k^2d}{gap}\right)$   ✗ | no | no |
| $k$-GenEV | LazyEV **Theorem 4.3** | $\widetilde{O}\left(\frac{knnz(B)\sqrt{\kappa_B}}{\sqrt{gap}} + \frac{knnz(A)+k^2d}{\sqrt{gap}}\right)$ | no | yes |
| | LazyEV **Theorem 4.4** | $\widetilde{O}\left(\frac{knnz(B)\sqrt{\kappa_B}}{\sqrt{\varepsilon}} + \frac{knnz(A)+k^2d}{\sqrt{\varepsilon}}\right)$ | yes | yes |

| Problem | Paper | Running time     (✗ for beaten) | gap-free? | stochastic? |
|---|---|---|---|---|
| | AppGrad [20] | $\widetilde{O}\left(\frac{knnz(X,Y)\cdot\kappa+k^2d}{gap}\right)$    (local conv.)   ✗ | no | no |
| | CCALin [15] | $\widetilde{O}\left(\frac{knnz(X,Y)\cdot\sqrt{\kappa}+k^2d}{gap}\right)$   ✗ | no | no |
| | CCALin [15] | $\widetilde{O}\left(\frac{knnz(X,Y)\cdot(1+\sqrt{\kappa'/n})+k^2d}{gap}\right)$   ✗ | no | yes |
| $k$-CCA | LazyCCA **Theorem 6.2** | $\widetilde{O}\left(\frac{knnz(X,Y)\cdot(1+\sqrt{\kappa'/n})+k^2d}{\sqrt{gap}}\right)$ | no | yes |
| | LazyCCA **Theorem 6.3** | $\widetilde{O}\left(\frac{knnz(X,Y)\cdot(1+\sqrt{\kappa'/n})+k^2d}{\sqrt{\varepsilon}}\right)$ | yes | yes |
| | LazyCCA **Theorem 6.2** | $\widetilde{O}\left(knnz(X,Y)\cdot\left(1+\frac{\sqrt{\kappa'}}{\sqrt{gap\cdot\sigma_k}\cdot(nnz(X,Y)/kd)^{1/4}}\right)\right)$ | no | doubly |
| | LazyCCA **Theorem 6.3** | $\widetilde{O}\left(knnz(X,Y)\cdot\left(1+\frac{\sqrt{\kappa'}}{\sqrt{\varepsilon\cdot\sigma_k}\cdot(nnz(X,Y)/kd)^{1/4}}\right)\right)$ | yes | doubly |

Table 2: Performance comparison on $k$-GenEV and $k$-CCA.
In GenEV, $gap = \frac{\lambda_k-\lambda_{k+1}}{\lambda_k} \in [0,1]$ and $\kappa_B = \frac{\lambda_{\max}(B)}{\lambda_{\min}(B)} > 1$.
In CCA, $gap = \frac{\sigma_k-\sigma_{k+1}}{\sigma_k} \in [0,1]$, $\kappa = \frac{\lambda_{\max}(\text{diag}\{S_{xx},S_{yy}\})}{\lambda_{\min}(\text{diag}\{S_{xx},S_{yy}\})} > 1$, $\kappa' = \frac{2\max_i\{\|X_i\|^2,\|Y_i\|^2\}}{\lambda_{\min}(\text{diag}\{S_{xx},S_{yy}\})} \in [\kappa, 2n\kappa]$, and $\sigma_k \in [0,1]$.

---

**Remark 1.** Stochastic methods depend on a modified condition number $\kappa'$. The reason $\kappa' \in [\kappa, 2n\kappa]$ is in Def. 2.4.
**Remark 2.** All non-stochastic CCA methods in this table have been outperformed because $1 + \sqrt{\kappa'/n} \le O(\kappa)$.
**Remark 3.** Doubly-stochastic methods are not necessarily interesting. We discuss them in Section 1.2.

---

## 1.2 Our Side Results on Doubly-Stochastic Methods

Recall that when considering acceleration, there are two parameters $\kappa$ and $1/gap$. One can also design stochastic methods with respect to both parameters $\kappa$ and $1/gap$, meaning that

$$\text{with a running time proportional to} \quad 1 + \frac{\sqrt{\kappa'/n^c}}{\sqrt{gap}} \quad \text{(constant } c \text{ is usually } 1/2)$$

instead of $\frac{1}{\sqrt{gap}} + \frac{\sqrt{\kappa'/n}}{\sqrt{gap}}$ (stochastic) or $\frac{\sqrt{\kappa}}{\sqrt{gap}}$ (non-stochastic). We call such methods *doubly-stochastic*.

Unfortunately, doubly-stochastic methods are usually slower than stochastic ones. Take 1-CCA as an example. The best stochastic running time (obtained exclusively by us) for 1-CCA is $nnz(X,Y) \cdot \widetilde{O}\left(\frac{1+\sqrt{\kappa'/n}}{\sqrt{gap}}\right)$. In contrast, if one uses a doubly-stochastic method —either [29] or our LazyCCA— the running time becomes $nnz(X,Y) \cdot \widetilde{O}\left(1 + \frac{\sqrt{\kappa'}/n^{1/4}}{\sqrt{gap}\cdot\sigma_1}\right)$. Therefore, for 1-CCA,

> doubly-stochastic methods are faster than stochastic ones **only when**   $\frac{\kappa'}{\sigma_1} \le o(n^{1/2})$  .

The above condition is usually *not* satisfied. For instance, (1) $\kappa'$ is usually around $n$ for most interesting data-sets, cf. the experiments of [25]; (2) $\kappa'$ is between $n^{1/2}$ and $100n$ in all the CCA experiments of [29]; and (3) by Def. 2.4 it satisfies $\kappa' \ge d$ so $\kappa'$ cannot be smaller than $o(n^{1/2})$ unless $d \ll n^{1/2}$.[4] Even worse, parameter $\sigma_1 \in [0,1]$ is usually much smaller than 1. Note that $\sigma_1$ is scaling invariant: even if one scales $X$ and $Y$ up by the same factor, $\sigma_1$ remains unchanged.

---

[4] Note that item (3) $\kappa' \ge d$ may not hold in the more general setting of CCA, see Remark 2.5.

**Nevertheless**, in order to compare our `LazyCCA` framework with *all* relevant prior works, we also obtain doubly-stochastic running times for $k$-CCA. Our running time matches that of [29] when $k = 1$, and no doubly-stochastic running time for $k > 1$ was known before our work.

### 1.3 Other Related Works

For the easier task of PCA and SVD, the first gap-free result was obtained by Musco and Musco [22] (or in the online setting by [3]), the first stochastic result was obtained by Shamir [26], and the first accelerated stochastic result was obtained by Garber *et al.* [13, 14]. The shift-and-invert preconditioning technique of Garber *et al.* is also used in this paper.

For another related problem PCR (principle component regression), we recently obtained an accelerated method [4] as opposed the previously non-accelerated one [12]; however, the acceleration techniques in [4] are not relevant to this paper.

For GenEV and CCA, many scalable algorithms have been designed recently [19–21, 28, 30]. However, as summarized by the authors of CCALin, these cited methods are more or less heuristics and do not have provable guarantees. Furthermore, for $k > 1$, the AppGrad result of [20] only provides local convergence guarantees and thus requires a warm-start whose computational complexity is not discussed in their paper.

Finally, our algorithms on GenEV and CCA are based on finding vectors one-by-one, which is advantageous in practice because one does not need $k$ to be known and can stop the algorithm whenever the eigenvalues (or correlation values) are too small. Known approaches for $k > 1$ cases (such as GenELin, CCALin, AppGrad) find all $k$ vectors at once, therefore requiring $k$ to be known beforehand. As a separate note, these known approaches do not need the user to know the desired accuracy a priori but our `LazyEV` and `LazyCCA` algorithms do.

## 2 Preliminaries

For a vector $x$ we denote by $\|x\|$ or $\|x\|_2$ the Euclidean norm of $x$. Given a matrix $A$ we denote by $\|A\|_2$ and $\|A\|_F$ respectively the spectral and Frobenius norms of $A$. For $q \geq 1$, we denote by $\|A\|_{S_q}$ the Schatten $q$-norm of $A$. We write $A \succeq B$ if $A, B$ are symmetric and $A - B$ is positive semi-definite (PSD), and write $A \succ B$ if $A, B$ are symmetric but $A - B$ is positive definite (PD). We denote by $\lambda_{\max}(M)$ and $\lambda_{\min}(M)$ the largest and smallest eigenvalue of a symmetric matrix $M$, and by $\kappa_M$ the condition number $\lambda_{\max}(M)/\lambda_{\min}(M)$ of a PSD matrix $M$.

Throughout this paper, for a matrix $A \in \mathbb{R}^{n \times d}$, we define $\mathsf{nnz}(A) \overset{\text{def}}{=} \max\{n, d, N\}$ where $N$ is the number of non-zero entries of $A$. For two matrices $X, Y$, we denote by $\mathsf{nnz}(X, Y) = \mathsf{nnz}(X) + \mathsf{nnz}(Y)$, and by $X_i$ or $Y_i$ the $i$-th row vector of $X$ or $Y$. We also use $\mathsf{poly}(x_1, x_2, \ldots, x_t)$ to represent a quantity that is asymptotically at most polynomial in terms of variables $x_1, \ldots, x_t$. Given a column orthonormal matrix $U \in \mathbb{R}^{n \times k}$, we denote by $U^\perp \in \mathbb{R}^{n \times (n-k)}$ the column orthonormal matrix consisting of an arbitrary basis in the space orthogonal to the span of $U$'s columns.

Given a PSD matrix $B$ and a vector $v$, the value $v^\top B v$ is the $B$-inner product. Two vectors $v, w$ satisfying $v^\top B w = 0$ are known as $B$-orthogonal. Given a PSD matrix $B$, we denote by $B^{-1}$ the Moore-Penrose pseudoinverse of $B$ which is also PSD, and denote by $B^{1/2}$ the matrix square root of $B$ (satisfying $B^{1/2} \succeq 0$). All occurrences of $B^{-1}$, $B^{1/2}$ and $B^{-1/2}$ are for analysis purpose only. When implementing our algorithms, it only requires one to multiply $B$ to a vector.

**Definition 2.1** (GenEV). *Given symmetric matrices $A, B \in \mathbb{R}^{d \times d}$ where $B$ is positive definite. The generalized eigenvectors of $A$ with respect to $B$ are $v_1, \ldots, v_d$, where each $v_i$ is*

$$v_i \in \arg\max_{v \in \mathbb{R}^d} \left\{ |v^\top A v| \quad \text{such that} \quad \left\{ \begin{array}{l} v^\top B v = 1 \\ v^\top B v_j = 0 \ \forall j \in [i-1] \end{array} \right. \right\}$$

*The corresponding generalized eigenvalues $\lambda_1, \ldots, \lambda_d$ satisfy $\lambda_i = v_i^\top A v_i$ which is possibly negative.*

Following the tradition of [13, 29], we assume *without loss of generality* that $\lambda_i \in [-1, 1]$.

**Definition 2.2** (CCA). *Given $X \in \mathbb{R}^{n \times d_x}, Y \in \mathbb{R}^{n \times d_y}$, letting $S_{xx} = \frac{1}{n} X^\top X$, $S_{xy} = \frac{1}{n} X^\top Y$, $S_{yy} = \frac{1}{n} Y^\top Y$, the canonical-correlation vectors are $\{(\phi_i, \psi_i)\}_{i=1}^r$ where $r = \min\{d_x, d_y\}$ and $\forall i$:*

$$(\phi_i, \psi_i) \in \arg\max_{\phi \in \mathbb{R}^{d_x}, \psi \in \mathbb{R}^{d_y}} \left\{ \phi^\top S_{xy} \psi \quad \text{such that} \quad \left\{ \begin{array}{l} \phi^\top S_{xx} \phi = 1 \ \wedge \ \phi^\top S_{xx} \phi_j = 0 \ \forall j \in [i-1] \\ \psi^\top S_{yy} \psi = 1 \ \wedge \ \psi^\top S_{yy} \psi_j = 0 \ \forall j \in [i-1] \end{array} \right. \right\}$$

*The corresponding canonical-correlation coefficients $\sigma_1, \ldots, \sigma_r$ satisfy $\sigma_i = \phi_i^\top S_{xy} \psi_i \in [0, 1]$.*

We emphasize here that $\sigma_i$ always lies in $[0, 1]$ and is scaling-invariant. When dealing with a CCA problem, we also denote by $d = d_x + d_y$.

**Lemma 2.3** (CCA to GenEV). *Given a CCA problem with matrices $X \in \mathbb{R}^{n \times d_x}, Y \in \mathbb{R}^{n \times d_y}$, and suppose the canonical-correlation vectors and coefficients are $\{(\phi_i, \psi_i, \sigma_i)\}_{i=1}^r$ where $r = \min\{d_x, d_y\}$. Define $A = \begin{pmatrix} 0 & S_{xy} \\ S_{xy}^\top & 0 \end{pmatrix}$ and $B = \begin{pmatrix} S_{xx} & 0 \\ 0 & S_{yy} \end{pmatrix}$. Then, the GenEV problem of $A$ with respect to $B$ has $2r$ eigenvalues $\{\pm \sigma_i\}_{i=1}^r$ and corresponding generalized eigenvectors $\left\{ \begin{pmatrix} \phi_i \\ \psi_i \end{pmatrix}, \begin{pmatrix} -\phi_i \\ \psi_i \end{pmatrix} \right\}_{i=1}^n$. The remaining $d_x + d_y - 2r$ eigenvalues are zeros.*

**Definition 2.4.** *In CCA, let $A$ and $B$ be as defined in Lemma 2.3. We define condition numbers*

$$\kappa \overset{\text{def}}{=} \frac{\lambda_{\max}(B)}{\lambda_{\min}(B)} \ \text{and} \ \kappa' \overset{\text{def}}{=} \frac{2 \max_i \{\|X_i\|^2, \|Y_i\|^2\}}{\lambda_{\min}(B)} \quad \text{(it must satisfy } \kappa' \in [\kappa, 2n\kappa] \text{ and } \kappa' \geq d.)$$

*Proof.* We have $\lambda_{\max}(B) \leq \text{Tr}(B) = \frac{1}{n} \sum_i \|X_i\|^2 + \|Y_i\|^2 \leq 2 \max_i \{\|X_i\|^2, \|Y_i\|^2\}$ and thus $\kappa' \geq \kappa$. Suppose without loss of generality that $\|X_1\|^2 = \max_i \{\|X_i\|^2, \|Y_i\|^2\}$; then, we have $\lambda_{\max}(B) \geq \frac{1}{n} \|X_1\|^2 = \frac{1}{n} \max_i \{\|X_i\|^2, \|Y_i\|^2\}$. This implies $\kappa' \leq 2n\kappa$. Finally, $2 \max_i \{\|X_i\|^2, \|Y_i\|^2\} \geq \text{Tr}(B) \geq d \lambda_{\min}(B)$ and therefore $\kappa' \geq d$. $\qquad \square$

**Remark 2.5.** We have followed the very original definition of CCA [15, 20] by letting $S_{xx} = \frac{1}{n} X^\top X$ and $S_{yy} = \frac{1}{n} Y^\top Y$. In contrast, one prior work [29] considered the slightly more general version $S_{xx} = \gamma_x I + \frac{1}{n} X^\top X$ and $S_{yy} = \gamma_y I + \frac{1}{n} Y^\top Y$ for some $\gamma_x, \gamma_y \geq 0$. All of the results in this paper continue to hold in this more general setting, but we refrain from doing so for notational simplicity. (The only difference is that the parameter $\kappa'$ will no longer satisfy $\kappa' \geq d$ in Def. 2.4.)

**Lemma 2.6.** *Given matrices $X \in \mathbb{R}^{n \times d_x}, Y \in \mathbb{R}^{n \times d_y}$, let $A$ and $B$ be as defined in Lemma 2.3. For every $w \in \mathbb{R}^d$, `Katyusha` method [1] finds a vector $w' \in \mathbb{R}^d$ satisfying $\|w' - B^{-1} A w\| \leq \varepsilon$*

$$\text{in time} \quad O\left( \text{nnz}(X, Y) \cdot \left(1 + \sqrt{\kappa'/n}\right) \cdot \log \frac{\kappa \|w\|^2}{\varepsilon} \right) \ .$$

**Algorithm 1** $\mathtt{AppxPCA}^{\pm}(\mathcal{A}, M, \delta_\times, \varepsilon, p)$

---

**Input:** $\mathcal{A}$, an approximate matrix inversion method; $M \in \mathbb{R}^{d \times d}$, a symmetric matrix satisfying $-I \preceq M \preceq I$; $\delta_\times \in (0, 0.5]$, a multiplicative error; $\varepsilon \in (0, 1)$, a numerical accuracy parameter; and $p \in (0, 1)$, the confidence parameter.

1:    $\widehat{w}_0 \leftarrow \mathtt{RanInit}(d)$; $s \leftarrow 0$; $\lambda^{(0)} \leftarrow 1 + \delta_\times$;         $\diamond$   $\widehat{w}_0$ *is a random unit vector, see Def. 3.2*

2:    $m_1 \leftarrow \left\lceil 4 \log\left(\frac{288 d\theta}{p^2}\right) \right\rceil$, $m_2 \leftarrow \left\lceil \log\left(\frac{36 d\theta}{p^2 \varepsilon}\right) \right\rceil$;      $\diamond$   $\theta$ *is the parameter of* $\mathtt{RanInit}$, *see Def. 3.2*

3:    $\widetilde{\varepsilon}_1 \leftarrow \frac{1}{64 m_1}\left(\frac{\delta_\times}{48}\right)^{m_1}$ and $\widetilde{\varepsilon}_2 \leftarrow \frac{\varepsilon}{8 m_2}\left(\frac{\delta_\times}{48}\right)^{m_2}$

4:    **repeat**           $\diamond$   *it satisfies* $m_1 = T^{\mathrm{PM}}(8, 1/32, p)$ *and* $m_2 = T^{\mathrm{PM}}(2, \varepsilon/4, p)$, *see Lemma B.1*

5:      $s \leftarrow s + 1$;

6:      **for** $t = 1$ **to** $m_1$ **do**

7:          Apply $\mathcal{A}$ to find $\widehat{w}_t$ satisfying $\left\| \widehat{w}_t - (\lambda^{(s-1)} I - M)^{-1} \widehat{w}_{t-1} \right\| \leq \widetilde{\varepsilon}_1$;

8:      $w_a \leftarrow \widehat{w}_{m_1} / \|\widehat{w}_{m_1}\|$;

9:      Apply $\mathcal{A}$ to find $v_a$ satisfying $\left\| v_a - (\lambda^{(s-1)} I - M)^{-1} w_a \right\| \leq \widetilde{\varepsilon}_1$;

10:     **for** $t = 1$ **to** $m_1$ **do**

11:        Apply $\mathcal{A}$ to find $\widehat{w}_t$ satisfying $\left\| \widehat{w}_t - (\lambda^{(s-1)} I + M)^{-1} \widehat{w}_{t-1} \right\| \leq \widetilde{\varepsilon}_1$;

12:     $w_b \leftarrow \widehat{w}_{m_1} / \|\widehat{w}_{m_1}\|$;

13:     Apply $\mathcal{A}$ to find $v_b$ satisfying $\left\| v_b - (\lambda^{(s-1)} I + M)^{-1} w_b \right\| \leq \widetilde{\varepsilon}_1$;

14:     $\Delta^{(s)} \leftarrow \frac{1}{2} \cdot \frac{1}{\max\{w_a^\top v_a, w_b^\top v_b\} - \widetilde{\varepsilon}_1}$ and $\lambda^{(s)} \leftarrow \lambda^{(s-1)} - \frac{\Delta^{(s)}}{2}$;

15:   **until** $\Delta^{(s)} \leq \frac{\delta_\times \lambda^{(s)}}{12}$

16:   $f \leftarrow s$;

17:   **if** the last $w_a^\top v_a \geq w_b^\top v_b$ **then**

18:     **for** $t = 1$ **to** $m_2$ **do**

19:        Apply $\mathcal{A}$ to find $\widehat{w}_t$ satisfying $\left\| \widehat{w}_t - (\lambda^{(f)} I - M)^{-1} \widehat{w}_{t-1} \right\| \leq \widetilde{\varepsilon}_2$;

20:     **return** $(+, w)$ where $w \overset{\text{def}}{=} \widehat{w}_{m_2} / \|\widehat{w}_{m_2}\|$.

21:   **else**

22:     **for** $t = 1$ **to** $m_2$ **do**

23:        Apply $\mathcal{A}$ to find $\widehat{w}_t$ satisfying $\left\| \widehat{w}_t - (\lambda^{(f)} I + M)^{-1} \widehat{w}_{t-1} \right\| \leq \widetilde{\varepsilon}_2$;

24:     **return** $(-, w)$ where $w \overset{\text{def}}{=} \widehat{w}_{m_2} / \|\widehat{w}_{m_2}\|$.

25:   **end if**

---

# 3   Leading Eigenvector via Two-Sided Shift-and-Invert

In this section we define $\mathtt{AppxPCA}^{\pm}$, the multiplicative approximation algorithm for computing the *two-sided* leading eigenvector of a symmetric matrix using the shift-and-invert preconditioning framework [13, 14]. Our pseudo-code Algorithm 1 is a modification of Algorithm 5 in [13].

The main differences between $\mathtt{AppxPCA}^{\pm}$ and Algorithm 5 of [13] are two-fold. First, given a symmetric matrix $M$, $\mathtt{AppxPCA}^{\pm}$ simultaneously considers an upper-bounding shift together with a lower-bounding shift, and try to invert both $\lambda I - M$ and $\lambda I + M$. This allows us to determine approximately how close $\lambda$ is to the largest *and* the smallest eigenvalues of $M$, and decrease $\lambda$ accordingly; in the end, it outputs an approximate eigenvector of $M$ that corresponds to a negative eigenvalue if needed. Second, we provide a multiplicative-error guarantee rather than additive as originally appeared in [13]. Without this multiplicative-error guarantee, our final running time will depend on $\frac{1}{\mathtt{gap} \cdot \lambda_{\max}(M)}$ rather than $\frac{1}{\mathtt{gap}}$.[5] Of course, we believe the bulk of the credit for conceiving

---

[5]This is why the SI method [29] depends on $\frac{1}{\mathtt{gap} \cdot \sigma_1}$ in Table 1.

`AppxPCA`$^{\pm}$ belongs to the original authors of [13, 14].

**Theorem 3.1** (`AppxPCA`$^{\pm}$)**.** *Let $M \in \mathbb{R}^{d \times d}$ be a symmetric matrix with eigenvalues $1 \geq \lambda_1 \geq \cdots \geq \lambda_d \geq -1$ and corresponding eigenvectors $u_1, \ldots, u_d$. Let $\lambda^* = \|M\|_2 = \max\{\lambda_1, -\lambda_d\}$. With probability at least $1 - p$, `AppxPCA`$^{\pm}$ produces a pair $(sgn, w)$ satisfying*

*if $sgn = +$, then* $\displaystyle\sum_{i \in [d], \lambda_i \leq (1 - \delta_\times/2)\lambda^*} (w^\top u_i)^2 \leq \varepsilon$ *and* $w^\top M w \geq (1 - \delta_\times/2)(1 - 3\varepsilon)\lambda^*$ *, and*

*if $sgn = -$, then* $\displaystyle\sum_{i \in [d], \lambda_i \geq -(1 - \delta_\times/2)\lambda^*} (w^\top u_i)^2 \leq \varepsilon$ *and* $w^\top M w \leq -(1 - \delta_\times/2)(1 - 3\varepsilon)\lambda^*$ *.*

*Furthermore, the total number of oracle calls to $\mathcal{A}$ is $O(\log(1/\delta_\times)m_1 + m_2)$, and each time we call $\mathcal{A}$ it satisfies that $\frac{\lambda_{\max}(\lambda^{(s)}I - M)}{\lambda_{\min}(\lambda^{(s)}I - M)}, \frac{\lambda_{\max}(\lambda^{(s)}I + M)}{\lambda_{\min}(\lambda^{(s)}I + M)} \in [1, \frac{96}{\delta_\times}]$ and $\frac{1}{\lambda_{\min}(\lambda^{(s)}I - M)}, \frac{1}{\lambda_{\min}(\lambda^{(s)}I + M)} \leq \frac{48}{\delta_\times \lambda^*}$.*

We remark here that, unlike the original shift-and-invert method which chooses a random (Gaussian) unit vector in Line 1 of `AppxPCA`$^{\pm}$, we allow this initial vector to be generated from an arbitrary $\theta$-conditioned random vector generator, defined as follows:

**Definition 3.2.** *An algorithm `RanInit`$(d)$ is a $\theta$-conditioned random vector generator if $w =$ `RanInit`$(d)$ is a $d$-dimensional unit vector and, for every $p \in (0, 1)$, every unit vector $u \in \mathbb{R}^d$, with probability at least $1 - p$, it satisfies $(u^\top w)^2 \leq \frac{p^2 \theta}{9d}$.*

This modification is needed in order to obtain our efficient implementations of GenEV and CCA algorithms. One can construct $\theta$-conditioned random vector generator as follows:

**Proposition 3.3.** *Given PSD matrix $B \in \mathbb{R}^{d \times d}$, if we set `RanInit`$(d) \stackrel{\text{def}}{=} \frac{B^{1/2}v}{(v^\top Bv)^{0.5}}$ where $v$ is a random Gaussian vector, then `RanInit`$(d)$ is a $\theta$-conditioned random vector generator for $\theta = \kappa_B$.*

*Proof of Proposition 3.3.* We have

$$(u^\top w)^2 = \frac{\text{Tr}(uu^\top vBv^\top)}{v^\top Bv} \overset{\text{①}}{\geq} \frac{\text{Tr}(uu^\top vBv^\top)}{\lambda_{\max}(B)} \overset{\text{②}}{\geq} \frac{\lambda_{\min}(B) \cdot \text{Tr}(uu^\top vv^\top)}{\lambda_{\max}(B)} = \theta(u^\top v)^2 \ .$$

Above, ① is because $v^\top Bv \leq \lambda_{\max}(B) \cdot \|v\|_2^2 = \lambda_{\max}(B)$, and ② follows from the fact that $vBv^\top \succeq v(\lambda_{\min}(B)I)v^\top = \lambda_{\min}(B)vv^\top$. Finally, using for instance [8, Lemma 5], it holds with probability at least $1 - p$ that $(u^\top v)^2 \geq \frac{p^2}{9d}$. $\qquad\square$

# 4 LazyEV: Our Algorithm for Generalized Eigendecomposition

In this section, we propose `LazyEV` (see Algorithm 2) to compute approximately the $k$ "leading" eigenvectors corresponding to the $k$ largest *absolute* eigenvalues of some symmetric matrix $M \in \mathbb{R}^{d \times d}$. Later, we shall solve the $k$-GenEV problem by setting $M = B^{-1/2}AB^{-1/2}$ and using `LazyEV` to find the top $k$ leading eigenvectors of $M$, which correspond to the top $k$ leading generalized eigenvectors of $A$ with respect to $B$.

Our algorithm `LazyEV` is formally stated in Algorithm 2. It applies $k$ times `AppxPCA`$^{\pm}$, each time with a multiplicative error $\delta_\times/2$, and projects the matrix $M$ into the orthogonal space with respect to the obtained leading eigenvector. We state our main approximation theorem below.

---

**Algorithm 2** $\mathtt{LazyEV}(\mathcal{A}, M, k, \delta_\times, \varepsilon_{\mathsf{pca}}, p)$

---

**Input:** $\mathcal{A}$, an approximate matrix inversion method; $M \in \mathbb{R}^{d \times d}$, a matrix satisfying $-I \preceq M \preceq I$; $k \in [d]$, the desired rank; $\delta_\times \in (0, 1)$, a multiplicative error; $\varepsilon_{\mathsf{pca}} \in (0, 1)$, a numerical accuracy parameter; and $p \in (0, 1)$, a confidence parameter.

1: $M_0 \leftarrow M; V_0 \leftarrow [];$
2: **for** $s = 1$ **to** $k$ **do**
3:      $v'_s \leftarrow \mathtt{AppxPCA}^{\pm}(\mathcal{A}, M_{s-1}, \delta_\times/2, \varepsilon_{\mathsf{pca}}, p/k);$
4:      $v_s \leftarrow \big((I - V_{s-1}V_{s-1}^\top)v'_s\big)/\big\|(I - V_{s-1}V_{s-1}^\top)v'_s\big\|;$           $\diamond$ *project $v'_s$ to $V_{s-1}^\perp$*
5:      $V_s \leftarrow [V_{s-1}, v_s];$
6:      $M_s \leftarrow (I - v_s v_s^\top)M_{s-1}(I - v_s v_s^\top)$           $\diamond$ *we also have $M_s = (I - V_s V_s^\top)M(I - V_s V_s^\top)$*
7: **end for**
8: **return** $V_k$.

---

**Theorem 4.1** (approximation of $\mathtt{LazyEV}$). *Let $M \in \mathbb{R}^{d \times d}$ be a symmetric matrix with eigenvalues $\lambda_1, \ldots, \lambda_d \in [-1, 1]$ and corresponding eigenvectors $u_1, \ldots, u_d$, and assume $|\lambda_1| \geq \cdots \geq |\lambda_d|$.*

*For every $k \in [d]$, $\delta_\times, p \in (0, 1)$, there exists some $\varepsilon_{\mathsf{pca}} \leq O\big(\mathsf{poly}(\delta_\times, \frac{|\lambda_1|}{|\lambda_{k+1}|}, \frac{1}{d})\big)$ such that[6] $\mathtt{LazyEV}$ outputs a (column) orthonormal matrix $V_k = (v_1, \ldots, v_k) \in \mathbb{R}^{d \times k}$ which, with probability at least $1 - p$, satisfies all of the following properties. (Denote by $M_s = (I - V_s V_s^\top)M(I - V_s V_s^\top)$.)*

(a) *Correlation guarantee: $\|V_k^\top U\|_2 \leq \varepsilon$,*
    *where $U = (u_j, \ldots, u_d)$ and $j$ is the smallest index satisfying $|\lambda_j| \leq (1 - \delta_\times)\|M_{k-1}\|_2$.*

(b) *Spectral norm guarantee: $|\lambda_{k+1}| \leq \|M_k\|_2 \leq \frac{|\lambda_{k+1}|}{1 - \delta_\times}$.*

(c) *Rayleigh quotient guarantee: $(1 - \delta_\times)|\lambda_k| \leq |v_k^\top M v_k| \leq \frac{1}{1 - \delta_\times}|\lambda_k|$.*

(d) *Schatten-$q$ norm guarantee: for every $q \geq 1$,*

$$\|M_k\|_{S_q} \leq \frac{(1 + \delta_\times)^2}{(1 - \delta_\times)^2}\Big(\sum_{i=k+1}^{d} \lambda_i^q\Big)^{1/q} = \frac{(1 + \delta_\times)^2}{(1 - \delta_\times)^2} \min_{V \in \mathbb{R}^{d \times k}, V^\top V = I} \big\{\|(I - VV^\top)M(I - VV^\top)\|_{S_q}\big\} \ .$$

The next theorem states that, if $M = B^{-1/2}AB^{-1/2}$, then $\mathtt{LazyEV}$ can be implemented without ever needing to compute $B^{1/2}$ or $B^{-1/2}$.

**Theorem 4.2** (running time of $\mathtt{LazyEV}$). *Let $A, B \in \mathbb{R}^{d \times d}$ be two symmetric matrices satisfying $B \succ 0$ and $-B \preceq A \preceq B$. Suppose $M = B^{-1/2}AB^{-1/2}$ and $\mathtt{RanInit}(d)$ is the random vector generator defined in Proposition 3.3 with respect to $B$. Then, the computation of $\mathbb{V} \leftarrow B^{-1/2}\mathtt{LazyEV}(\mathcal{A}, M, k, \delta_\times, \varepsilon_{\mathsf{pca}}, p)$ can be implemented to run in time*

- $\widetilde{O}\Big(\frac{k\mathsf{nnz}(B) + k^2 d + k\Upsilon}{\sqrt{\delta_\times}}\Big)$ *where $\Upsilon$ is the time to multiply $B^{-1}A$ to a vector,[7] or*

- $\widetilde{O}\Big(\frac{k\sqrt{\kappa_B}\mathsf{nnz}(B) + k\mathsf{nnz}(A) + k^2 d}{\sqrt{\delta_\times}}\Big)$ *if we use Conjugate gradient to multiply $B^{-1}A$ to a vector.*

*Above, the $\widetilde{O}$ notation hides polylogarithmic factors with respect to $1/\varepsilon_{\mathsf{pca}}, 1/\delta_\times, 1/p, \kappa_B, d$.*

Our main theorems immediately imply the following corollaries (proved in Appendix E.2):

---

[6]The complete specifications of $\varepsilon_{\mathsf{pca}}$ is included in Appendix E. Since our final running time only depends on $\log(1/\varepsilon_{\mathsf{pca}})$, we have not attempted to improve the constants in this polynomial dependency.

[7]More precisely, to compute $(B^{-1}A)w$ for some vector $w$ with error $\varepsilon$ where $\log(1/\varepsilon) = \widetilde{O}(1)$.

**Theorem 4.3** (gap-dependent $k$-GenEV). *Let $A, B \in \mathbb{R}^{d \times d}$ be two symmetric matrices satisfying $B \succ 0$ and $-B \preceq A \preceq B$. Suppose the generalized eigenvalue and eigenvector pairs of $A$ with respect to $B$ are $\{(\lambda_i, u_i)\}_{i=1}^{d}$, and it satisfies $1 \geq |\lambda_1| \geq \cdots \geq |\lambda_d|$. Let $\mathsf{gap} = \frac{|\lambda_k| - |\lambda_{k+1}|}{|\lambda_k|} \in [0, 1]$ be the relative eigengap. For fixed $\varepsilon, p > 0$, consider the output*

$$\mathbb{V}_k \leftarrow B^{-1/2} \mathtt{LazyEV}\left( \mathcal{A}, B^{-1/2}AB^{-1/2}, k, \mathsf{gap}, O\left( \frac{\varepsilon^4 \cdot \mathsf{gap}}{k^3 (\sigma_1/\sigma_k)^4} \right), p \right) \in \mathbb{R}^{d \times k} \ .$$

*Then, defining $\mathbb{W} = (u_{k+1}, \ldots, u_d)$, we have with probability at least $1 - p$:*

$$\mathbb{V}_k^{\top} B \mathbb{V}_k = I \quad and \quad \|\mathbb{V}_k^{\top} B \mathbb{W}\|_2 \leq \varepsilon \ .$$

*Moreover, our running time is $\widetilde{O}\left( \frac{k\sqrt{\kappa_B}\mathsf{nnz}(B) + k\mathsf{nnz}(A) + k^2 d}{\sqrt{\mathsf{gap}}} \right)$.*

---

**Theorem 4.4** (gap-free $k$-GenEV). *In the same setting as Theorem 4.3, for $\varepsilon, p > 0$, consider*

$$(\mathtt{v}_1, \ldots, \mathtt{v}_k) \overset{\text{def}}{=} \mathbb{V}_k \leftarrow B^{-1/2} \mathtt{LazyEV}\left( \mathcal{A}, B^{-1/2}AB^{-1/2}, k, \varepsilon, O\left( \frac{\varepsilon^5}{k^3 d^4 (\sigma_1/\sigma_{k+1})^{12}} \right), p \right) \ .$$

*Then, with probability at least $1 - p$: it satisfies $\mathbb{V}_k^{\top} B \mathbb{V}_k = I$ and*

$$\forall s \in [k] \colon |\mathtt{v}_s^{\top} A \mathtt{v}_s| \in \left[ (1 - \varepsilon)|\lambda_s|, \frac{|\lambda_s|}{1 - \varepsilon} \right] \quad and \quad \max_{\mathtt{w} \in \mathbb{R}^d \wedge \mathtt{w}^{\top} B \mathbb{V}_k = 0} \left| \frac{\mathtt{w}^{\top} A \mathtt{w}}{\mathtt{w}^{\top} B \mathtt{w}} \right| \leq \frac{|\lambda_{k+1}|}{1 - \varepsilon} \ .$$

*Moreover, our running time is $\widetilde{O}\left( \frac{k\sqrt{\kappa_B}\mathsf{nnz}(B) + k\mathsf{nnz}(A) + k^2 d}{\sqrt{\varepsilon}} \right)$.*

---

# 5 Ideas Behind Theorems 4.1 and 4.2

Our $\mathtt{LazyEV}$ algorithm reduces the problem of finding generalized eigenvectors to finding regular eigenvectors of $M = B^{-1/2}AB^{-1/2}$. In Section 5.1 we discuss how to ensure accuracy: that is, why does $\mathtt{LazyEV}$ guarantee to find approximately the top *absolute* eigenvectors of $M$; and in Section 5.2 we discuss how to implement $\mathtt{LazyEV}$ without ever needing to compute $B^{1/2}$ or $B^{-1/2}$.

## 5.1 Ideas Behind Theorem 4.1: Approximation Guarantee of GenEV

Our approximation guarantee in Theorem 4.1 is a natural generalization of the recent work on fast iterative methods to find the top $k$ eigenvectors of a PSD matrix $M$ [2]. That method is called $\mathtt{LazySVD}$. At a high level, $\mathtt{LazySVD}$ finds the top $k$ eigenvectors of $M$ one by one but only *approximately*. Starting with $M_0 = M$, in the $s$-th iteration where $s \in [k]$, $\mathtt{LazySVD}$ computes approximately the leading eigenvector of matrix $M_{s-1}$ (using shift-and-invert [13]) and call it $v_s$. Then, $\mathtt{LazySVD}$ projects $M_s \leftarrow (I - v_s v_s^{\top}) M_{s-1} (I - v_s v_s^{\top})$ and proceeds to the next iteration.

While the algorithmic idea of $\mathtt{LazySVD}$ is simple, the analysis requires some careful linear algebraic lemmas. Most notably, if $v_s$ is an approximate leading eigenvector of $M_{s-1}$, then one needs to prove that the small eigenvectors of $M_{s-1}$ somehow still "embed" into that of $M_s$ after projection. This is achieved by a gap-free variant of the Wedin theorem plus a few other technical lemmas, and we recommend interested readers to see the high-level overview section of [2].

In this paper, to relax the assumption that $M$ is PSD, and to find leading eigenvectors whose *absolute* eigenvalues are large, we have to make some non-trivial changes in the algorithm and the analysis. On the algorithm side, $\mathtt{LazyEV}$ replaces the use of the shift-and-invert protocol in $\mathtt{LazySVD}$ with our two-sided variant developed in Section 3. On the analysis side, we have to make sure all lemmas properly deal with negative eigenvalues: for instance, if we perform a projection $M' \leftarrow (I - vv^{\top}) M (I - vv^{\top})$ where $v$ correlates by at most $\varepsilon$ with all eigenvectors of $M$ whose absolute eigenvalues are smaller than a threshold $\mu$, then, after the projection, we need to prove that these eigenvectors can be approximately "embedded" into the eigenspace spanned by all eigenvectors of $M'$ whose *absolute* eigenvalues are smaller than $\mu + \tau$. The approximation of this embedding should depend on $\varepsilon, \mu$ and $\tau$. See Lemma C.4 in the appendix.

The proof of Theorem 4.1 is in Appendix E, and the matrix algebraic lemmas are in Appendix C.

## 5.2 Proof of Theorem 4.2: Fast Implementation of GenEV

We can implement `LazyEV` efficiently without the necessity of computing $B^{1/2}$ or $B^{-1/2}$. In each iteration of `LazyEV`, we call `AppxPCA`$^{\pm}$ and compute a vector $v'_s$. We do not explicitly store $v'_s$, but rather write it as $v'_s = B^{1/2}\mathrm{v}'_s$ and store only $\mathrm{v}'_s \in \mathbb{R}^d$. We shall later ensure that `AppxPCA`$^{\pm}$ outputs $v'_s$ directly. Similarly, we also write $v_s = B^{1/2}\mathrm{v}_s$ and only store $\mathrm{v}_s$. All together, we do not explicitly compute $V_s$, but instead write $V_s = B^{1/2}\mathbb{V}_s$ and only keep track of $\mathbb{V}_s \in \mathbb{R}^{d \times s}$.

Now, the computation of $v_s$ becomes the $B$-projection into the $\mathbb{V}_{s-1}$ space:

$$\|(I - V_{s-1}V_{s-1}^{\top})v'_s\| = \|B^{1/2}\mathrm{v}'_s - B^{1/2}\mathbb{V}_{s-1}\mathbb{V}_{s-1}^{\top}B\mathrm{v}'_s\| = \left(\left(\mathrm{v}'_s - \mathbb{V}_{s-1}\mathbb{V}_{s-1}^{\top}B\mathrm{v}'_s\right)^{\top}B\left(\mathrm{v}'_s - \mathbb{V}_{s-1}\mathbb{V}_{s-1}^{\top}B\mathrm{v}'_s\right)\right)^{1/2}$$

$$\|(I - V_{s-1}V_{s-1}^{\top})v'_s\| \cdot \mathrm{v}_s = B^{-1/2}(I - V_{s-1}V_{s-1}^{\top})v'_s = B^{-1/2}(I - V_{s-1}V_{s-1}^{\top})B^{1/2}\mathrm{v}'_s = \mathrm{v}'_s - \mathbb{V}_{s-1}\mathbb{V}_{s-1}^{\top}B\mathrm{v}'_s$$

and this can be implemented to run in $O(kd + \mathtt{nnz}(B))$ time. Finally, we write

$$M_s = (I - V_sV_s^{\top})B^{-1/2}AB^{-1/2}(I - V_sV_s^{\top}) = B^{-1/2}(I - B\mathbb{V}_s\mathbb{V}_s^{\top})A(I - \mathbb{V}_s\mathbb{V}_s^{\top}B)B^{-1/2}$$

and only pass it implicity to `AppxPCA`$^{\pm}$ (without directly computing this matrix).

To implement `AppxPCA`$^{\pm}$, we again write all vectors $\widehat{w}_t = B^{1/2}\mathrm{w}_t$ and only store $\mathrm{w}_t$. Thus, the normalization $w_a \leftarrow \widehat{w}_{m_1}/\|\widehat{w}_{m_1}\|_2$ becomes the $B$-normalization $\mathrm{w}_a \leftarrow \mathrm{w}_{m_1}/(\mathrm{w}_{m_1}^{\top}B\mathrm{w}_{m_1})^{1/2}$ which runs in $O(\mathtt{nnz}(B))$ time. Recall that `AppxPCA`$^{\pm}$ makes a polylogarithmic number of calls to the matrix inversion subroutine $\mathcal{A}$, each time requesting to approximately invert either $\lambda I - M_s$ or $\lambda I + M_s$. Let us only focus on inverting $\lambda I - M_s$ and the other case is similar. We write

$$\mathbb{N} \overset{\text{def}}{=} B^{-1/2}(\lambda I - M_s)B^{1/2} = \lambda I - (I - \mathbb{V}_s\mathbb{V}_s^{\top}B)B^{-1}A(I - \mathbb{V}_s\mathbb{V}_s^{\top}B) \ .$$

Now, the accuracy requirement in `AppxPCA`$^{\pm}$ becomes

$$\text{find } \mathrm{w}_t \text{ satisfying } \|B^{1/2}\mathrm{w}_t - (\lambda I - M_s)^{-1}B^{1/2}\mathrm{w}_{t-1}\| \leq \widetilde{\varepsilon}$$
$$\Longleftrightarrow \text{ find } \mathrm{w}_t \text{ satisfying } \|B^{1/2}\mathrm{w}_t - B^{1/2}\mathbb{N}^{-1}\mathrm{w}_{t-1}\| \leq \widetilde{\varepsilon}$$
$$\Longleftarrow \text{ find } \mathrm{w}_t \text{ satisfying } \|\mathrm{w}_t - \mathbb{N}^{-1}\mathrm{w}_{t-1}\| \leq \widetilde{\varepsilon}/\sqrt{\lambda_{\max}(B)}$$

Using accelerated gradient descent (see Theorem D.1), we can reduce this approximate inversion $\mathrm{w}_t \leftarrow \mathbb{N}^{-1}\mathrm{w}_{t-1}$ to $T$ times of approximate matrix-vector multiplication (i.e., $w' \leftarrow \mathbb{N}w$) for $T = \widetilde{O}(\sqrt{\kappa_{(\lambda I - M_s)}})$.[8] We can further derive that $T = \widetilde{O}(1/\sqrt{\delta_{\times}})$ owing to Theorem 3.1. Notice that each time we compute $w' \leftarrow \mathbb{N}w$ it suffices to compute it to an additive accuracy $\|w' - \mathbb{N}w\| \leq \varepsilon$ where the error satisfies $\log(1/\varepsilon) = \widetilde{O}(1)$.

Finally, the matrix-vector multiplication $\mathbb{N}w = \lambda w - (I - \mathbb{V}_s\mathbb{V}_s^{\top}B)B^{-1}A(I - \mathbb{V}_s\mathbb{V}_s^{\top}B)w$ consists of two rank-$s$ $B$-projections which run in time $O(\mathtt{nnz}(B) + kd)$, plus the time needed to multiply $B^{-1}A$ to a vector. This proves that `LazyEV` can be implemented so that

- It computes matrix-vector multiplication of the form $w' \leftarrow B^{-1}Aw$ a total of $\widetilde{O}(k/\sqrt{\delta_{\times}})$ times, each time to an accuracy $\varepsilon$ where $\log(1/\varepsilon) = \widetilde{O}(1)$;
- The rest of the computation costs a total of $\widetilde{O}\left((k\mathtt{nnz}(B) + k^2d)/\sqrt{\delta_{\times}}\right)$ time.

This finishes the proof of the first item of the theorem.

---

[8] This reduction would be obvious if we required the matrix-vector multiplication to be exact, and for instance Chebyshev method serves for exactly this purpose. However, in order to relax the multiplication to be approximate, and without incurring an error that blows up exponentially with $T$, we build our own inexact variant of the accelerate gradient descent method `AGD`$^{\text{inexact}}$ in Appendix D that could be of independent interest.

---

**Algorithm 3** `LazyCCA`$(\mathcal{A}, M, k, \delta_\times, \varepsilon_{\mathsf{pca}}, p)$

---

**Input:** $\mathcal{A}$, an approximate matrix inversion method;
  $M \in \mathbb{R}^{d \times d}$, a matrix satisfying $-I \preceq M \preceq I$;
  $k \in [d]$, the desired rank;
  $\delta_\times \in (0, 1)$, a multiplicative error;
  $\varepsilon_{\mathsf{pca}} \in (0, 1)$, a numerical accuracy parameter; and
  $p \in (0, 1)$, a confidence parameter.

1: $M_0 \leftarrow M$;
2: $V_0 = [\,]$;
3: **for** $s = 1$ to $k$ **do**
4:   $v'_s \leftarrow \mathtt{AppxPCA}^\pm(\mathcal{A}, M_{s-1}, \delta_\times/2, \varepsilon_{\mathsf{pca}}, p/k)$;
5:   $v_s \leftarrow \big((I - V_{s-1} V_{s-1}^\top) v'_s\big)/\big\|(I - V_{s-1} V_{s-1}^\top) v'_s\big\|$;   $\diamond$ *project $v'_s$ to $V_{s-1}^\perp$*
6:   write $v_s = (\xi'_s, \zeta'_s)$ where $\xi'_s \in \mathbb{R}^{d_x}$ and $\zeta'_s \in \mathbb{R}^{d_y}$;
7:   $\xi_s \leftarrow \xi'_s/(\sqrt{2}\|\xi'_s\|_2)$ and $\zeta_s \leftarrow \zeta'_s/(\sqrt{2}\|\zeta'_s\|_2)$;
8:   $V_s \leftarrow \left[V_{s-1}, \begin{pmatrix} \xi_s & -\xi_s \\ \zeta_s & \zeta_s \end{pmatrix}\right]$;
9:   $M_s \leftarrow \big(I - 2\mathrm{diag}(\xi_s \xi_s^\top, \zeta_s \zeta_s^\top)\big) M_{s-1} \big(I - 2\mathrm{diag}(\xi_s \xi_s^\top, \zeta_s \zeta_s^\top)\big)$
     $\diamond$ *or equivalently, $M_s = (I - V_s V_s^\top) M (I - V_s V_s^\top)$*
10: **end for**
11: **return** $V_k$.

---

As for the second item, we simply notice that whenever we want to compute $w' \leftarrow B^{-1} A w$, we can first compute $Aw$ in time $O(\mathtt{nnz}(A))$, and then use Conjugate gradient [27] to compute $B^{-1}$ applied to this vector. The running time of Conjugate gradient is at most $\widetilde{O}\big(\sqrt{\kappa_B} \cdot \mathtt{nnz}(B)\big)$ where the $\widetilde{O}$ factor hides a logarithmic factor on the accuracy. □

# 6  LazyCCA: Our Algorithm for Canonical Correlation Analysis

We propose `LazyCCA` (see Algorithm 3), a variant of `LazyEV` that is specially designed for matrices $M$ of the form $M = B^{-1/2} A B^{-1/2}$, where $A$ and $B$ come from a CCA instance following Lemma 2.3.

More specifically, recall from Lemma 2.3 that the eigenvectors of matrices $M$ arising from CCA instances are symmetric: if $(\xi, \zeta)$ is a normalized eigenvector of $M$ with eigenvalue $\sigma$ where $\xi \in \mathbb{R}^{d_x}, \zeta \in \mathbb{R}^{d_y}$, then $(-\xi, \zeta)$ is also a normalized eigenvector but with eigenvalue $-\sigma$. Furthermore, since $(\xi, \zeta)$ is orthogonal to $(-\xi, \zeta)$, we must have $\|\xi\| = \|\zeta\| = 1/\sqrt{2}$. Our `LazyCCA` method is designed to ensure such symmetry and orthogonality as well. When an approximate eigenvector $v_s = (\xi'_s, \zeta'_s)$ is obtained, we re-scale the pair to $\xi_s$ and $\zeta_s$ where both of them have norm exactly $1/\sqrt{2}$ (see Line 7 of `LazyCCA`). Then, we simultaneously add two (orthogonal) approximate eigenvectors $(\xi_s, \zeta_s)$ and $(-\xi_s, \zeta_s)$ to the column orthonormal matrix $V_s$.

**Remark 6.1.** This re-scaling step, together with the fact that we find vector pairs one by one, allows us to provide *per-vector* guarantee on the obtained approximate correlation vectors (see Theorem 6.3). This is in contrast to CCALin which is based on subspace power method so can only find the subspace spanned by the top $k$ correlation vectors but not distinguish them.

We prove three main theorems in the appendix:

- A main theorem for the approximation guarantee (see Theorem F.1).

  The main "delta" between the proofs of Theorem F.1 and Theorem 4.1 is to show that, after re-scaling, the vector $(\xi_s, \zeta_s)$ is also an approximate leading eigenvector of $M_{s-1}$. In particular, its Rayleigh quotient can only become better after scaling (see (F.3) in the appendix).

- A main theorem for the stochastic running time (see Theorem F.2).

  The main "delta" between the proofs of Theorem F.2 and Theorem 4.2 is to replace the use of conjugate gradient with the `Katyusha` method to multiply $B^{-1}A$ to a vector (see Lemma 2.6).

- A main theorem for the doubly-stochastic running time guarantee (see Theorem F.3).

  The main "delta" between the proofs of Theorem F.3 and Theorem F.2 is to use accelerated SVRG as opposed to `Katyusha` to compute matrix inverse.

We state below the final statements on `LazyCCA`.

---

**Theorem 6.2** (gap-dependent $k$-CCA). *Let $X \in \mathbb{R}^{n \times d_x}, Y \in \mathbb{R}^{n \times d_y}$ be two matrices with canonical-correlation coefficients $1 \geq \sigma_1 \geq \cdots \sigma_r \geq 0$ and the corresponding correlation vectors $\{(\phi_i, \psi_i)\}_{i=1}^r$. Let $\mathsf{gap} = \frac{\sigma_k - \sigma_{k+1}}{\sigma_k} \in [0, 1]$ be the relative gap, and define $A = [[0, S_{xy}]; [S_{xy}^\top, 0]]$ and $B = \mathrm{diag}(S_{xx}, S_{yy})$ following Def. 2.2. For every $\varepsilon, p > 0$, consider the output*

$$\begin{pmatrix} \pm\phi_1' & \dots & \pm\phi_k' \\ \psi_1' & \dots & \psi_k' \end{pmatrix} \stackrel{\text{def}}{=} \mathbb{V}_k \leftarrow \sqrt{2} B^{-1/2} \texttt{LazyCCA}\left(\mathcal{A}, B^{-1/2}AB^{-1/2}, \varepsilon, \mathsf{gap}, O\big(\tfrac{\varepsilon^4 \cdot \mathsf{gap}}{k^3(\sigma_1/\sigma_k)^4}\big), p\right) \ .$$

*Then, letting*

$$\mathbb{V}_\phi = (\phi_1', \dots, \phi_k'), \ \ \mathbb{V}_\psi = (\psi_1', \dots, \psi_k'), \ \ \mathbb{W}_\phi = (\phi_{k+1}, \phi_{k+2}, \dots), \ \ and \ \ \mathbb{W}_\psi = (\psi_{k+1}, \psi_{k+2}, \dots) \ ,$$

*we have with probability at least $1 - p$:*

$$\mathbb{V}_\phi \in \mathbb{R}^{d_x \times k} \ satisfies \ \mathbb{V}_\phi^\top S_{xx} \mathbb{V}_\phi = I \ and \ \|\mathbb{V}_\phi^\top S_{xx} \mathbb{W}_\phi\|_2 \leq \varepsilon \ ,$$

$$\mathbb{V}_\psi \in \mathbb{R}^{d_y \times k} \ satisfies \ \mathbb{V}_\psi^\top S_{yy} \mathbb{V}_\psi = I \ and \ \|\mathbb{V}_\psi^\top S_{yy} \mathbb{W}_\psi\|_2 \leq \varepsilon \ .$$

*The (stochastic) running time is $\widetilde{O}\big(\frac{\mathsf{nnz}(X,Y) \cdot (1 + \sqrt{\kappa'/n}) + k^2 d}{\sqrt{\mathsf{gap}}}\big)$. The doubly-stochastic running time is $\widetilde{O}\big(\mathsf{nnz}(X, Y) \cdot \big(1 + \frac{\sqrt{\kappa'/n^{1/4}}}{\sqrt{\mathsf{gap} \cdot \sigma_1}}\big)\big)$ for 1-CCA, and $\widetilde{O}\big(k\mathsf{nnz}(X, Y) \cdot \big(1 + \frac{\sqrt{\kappa'}}{\sqrt{\mathsf{gap} \cdot \sigma_k} \cdot (\mathsf{nnz}(X,Y)/kd)^{1/4}}\big)\big)$ for $k$-CCA as long as $\mathsf{nnz}(X, Y) \geq kd$.*

---

**Theorem 6.3** (gap-free $k$-CCA). *In the same setting as Theorem 6.2, for $\varepsilon, p > 0$, consider*

$$\begin{pmatrix} \pm\phi_1' & \dots & \pm\phi_k' \\ \psi_1' & \dots & \psi_k' \end{pmatrix} = \mathbb{V}_k \leftarrow \sqrt{2} B^{-1/2} \texttt{LazyCCA}\left(\mathcal{A}, B^{-1/2}AB^{-1/2}, \varepsilon, \mathsf{gap}, O\big(\tfrac{\varepsilon^4 \cdot \mathsf{gap}}{k^3(\sigma_1/\sigma_k)^4}\big), p\right) \ .$$

*Letting $\mathbb{V}_\phi = (\phi_1', \dots, \phi_k') \in \mathbb{R}^{d_x \times k}$ and $\mathbb{V}_\psi = (\psi_1', \dots, \psi_k') \in \mathbb{R}^{d_y \times k}$, with probability at least $1 - p$,*

- $\mathbb{V}_\phi^\top S_{xx} \mathbb{V}_\phi = I$, $\mathbb{V}_\psi^\top S_{yy} \mathbb{V}_\psi = I$;
- $(1 - \varepsilon)\sigma_i \leq |\phi_i' S_{xy} \psi_i| \leq (1 + \varepsilon)\sigma_i$ for every $i \in [k]$; and
- $\max_{\phi \in \mathbb{R}^{d_x}, \psi \in \mathbb{R}^{d_y}} \left\{ \phi^\top S_{xy} \psi \ \middle| \ \phi^\top S_{xx} \mathbb{V}_\phi = 0 \wedge \psi^\top S_{yy} \mathbb{V}_\psi = 0 \right\} \leq (1 + \varepsilon)\sigma_{k+1} \ .$

*The (stochastic) running time is $\widetilde{O}\big(\frac{\mathsf{nnz}(X,Y) \cdot (1 + \sqrt{\kappa'/n}) + k^2 d}{\sqrt{\varepsilon}}\big)$. The doubly-stochastic running time is $\widetilde{O}\big(\mathsf{nnz}(X, Y) \cdot \big(1 + \frac{\sqrt{\kappa'/n^{1/4}}}{\sqrt{\varepsilon} \cdot \sigma_1}\big)\big)$ for 1-CCA, and $\widetilde{O}\big(k\mathsf{nnz}(X, Y) \cdot \big(1 + \frac{\sqrt{\kappa'}}{\sqrt{\varepsilon} \cdot \sigma_k \cdot (\mathsf{nnz}(X,Y)/kd)^{1/4}}\big)\big)$ for $k$-CCA as long as $\mathsf{nnz}(X, Y) \geq kd$.*

---

**Remark 6.4.** The doubly-stochastic running times for the general $k > 1$ setting can be slightly faster (but notationally more involved) than the ones stated above. For instance, we are aware of a proof that gives running time

$$\widetilde{O}\big(k\mathsf{nnz}(X, Y) \cdot \big(1 + \frac{\sqrt{\kappa'}}{\sqrt{\mathsf{gap} \cdot \sigma_k}} \cdot \frac{\max\{k \min\{\sigma_1, \sigma_1 \kappa/d\}, 1\}^{1/4} d^{1/4}}{(\mathsf{nnz}(X, Y))^{1/4}}\big)\big)$$

for the gap-dependent case of $k$-CCA, or the same formula but replacing gap with $\varepsilon$ for the gap-free case. We refrain from proving it in full because as discussed in Section 1.2, doubly-stochastic running times are not necessarily interesting.

# Appendix

## A  Proof of Lemma 2.6

*Proof of Lemma 2.6.* First of all, computing $B^{-1}Aw$ is equivalent to minimizing $f(x) \overset{\text{def}}{=} \frac{1}{2}x^\top Bx - x^\top Aw$. Suppose we write $x = (x_1, x_2)$ and $w = (w_1, w_2)$ where $x_1, w_1 \in \mathbb{R}^{d_x}, x_2, w_2 \in \mathbb{R}^{d_y}$, then one can rewrite $f$ as $f(x) = \frac{1}{2n}\left(\|Xx_1 - Yw_2\|_2^2 + \|Xw_1 - Yx_2\|_2^2\right) + C$ where $C$ is a fixed constant. Therefore, one can also write

$$f(x) = \frac{1}{2n}\sum_{i=1}^{n}(\langle X_i, x_1\rangle - \langle Y_i, w_2\rangle)^2 + (\langle Y_i, x_2\rangle - \langle X_i, w_1\rangle)^2$$

where each $X_i$ is a row vector of $X$ and $Y_i$ is a row vector of $Y$. In such a case, we observe that

- $f(x)$ is an average of $2n$ smooth functions, where function $(\langle X_i, x_1\rangle - \langle Y_i, w_2\rangle)^2$ is smooth with parameter $2\|X_i\|^2$ (meaning Hessian bounded by $\|X_i\|^2$ in spectral norm) and $(\langle Y_i, x_2\rangle - \langle X_i, w_1\rangle)^2$ is smooth with parameter $2\|Y_i\|^2$. In other words, each function has a smoothness parameter at most $2\max_i\{\|X_i\|^2, \|Y_i\|^2\}$.

- $f(x)$ is at least $\lambda_{\min}(B) = \min\{\lambda_{\min}(S_{xx}), \lambda_{\min}(S_{yy})\}$ strongly convex, meaning the Hessian $\nabla^2 f(x)$ has no eigenvalue less than this quantity.

For such reason, one can apply the convergence theorem of Katyusha [1] to find some $x$ such that $f(x) - \min_y f(y) \leq \widetilde{\varepsilon}$ in time $O\left(\mathsf{nnz}(X, Y) \cdot \left(1 + \sqrt{\kappa'/n}\right)\log\frac{f(x^0) - f(x^*)}{\widetilde{\varepsilon}}\right)$ where $x^0$ is an arbitrary starting vector fed into Katyusha and $x^*$ is the exact minimizer. If we choose $x^0$ to be the zero vector, it is easy to verify that $f(x^0) - f(x^*) \leq O(\lambda_{\max}(B) \cdot \|w\|^2)$.

It is not hard to see that an additive $\widetilde{\varepsilon}$ minimizer of $f(x)$ implies an $\varepsilon$-approximate solution for the inverse $\|w' - B^{-1}Aw\| \leq \varepsilon$ where $\varepsilon^2 = 2\widetilde{\varepsilon}/\lambda_{\min}(B)$. This finishes the proof of the running time in Lemma 2.6. $\qquad\square$

## B  Proof for Section 3: Two-Sided Shift-and-Invert

### B.1  Inexact Power Method

In this subsection we review some classical convergence lemmas regarding power method and its inexact variant. These lemmas almost directly follow from previous results such as [13, 15], and are more similar to [2]. We skip the proofs in this paper.

Consider power method that starts with a random unit vector $w_0 \leftarrow \mathtt{RanInit}(d)$ and apply $w_t \leftarrow Mw_{t-1}/\|Mw_{t-1}\|$ iteratively.

**Lemma B.1** (Exact Power Method). *Let $M$ be a PSD matrix with eigenvalues $\lambda_1 \geq \cdots \geq \lambda_d$ and the corresponding eigenvectors $u_1, \ldots, u_d$. Fix an error tolerance $\varepsilon > 0$, parameter $\kappa \geq 1$, and failure probability $p > 0$, define*

$$T^{\mathrm{PM}}(\kappa, \varepsilon, p) = \left\lceil \frac{\kappa}{2}\log\left(\frac{9d\theta}{p^2\varepsilon}\right)\right\rceil$$

*Then, with probability at least $1 - p$ it holds that $\forall t \geq T^{\mathrm{PM}}(\kappa, \varepsilon, p)$:*

$$\sum_{i \in [d], \lambda_i \leq (1 - 1/\kappa)\lambda_1} (w^\top u_i)^2 \leq \varepsilon \quad and \quad w^\top M w \geq (1 - 1/\kappa - \varepsilon)\lambda_1 \ .$$

**Lemma B.2** (Lemma 4.1 of [13])**.** *Let $M$ be a PSD matrix with eigenvalues $\lambda_1 \geq \cdots \lambda_d$. Fix an accuracy parameter $\widetilde{\varepsilon} > 0$, and consider two update sequences*

$$\widehat{w}_0^* = w_0, \qquad \forall t \geq 1 \colon \widehat{w}_t^* \leftarrow M \widehat{w}_{t-1}^*$$
$$\widehat{w}_0 = w_0, \qquad \forall t \geq 1 \colon \widehat{w}_t \text{ satisfies } \|\widehat{w}_t - M\widehat{w}_{t-1}\| \leq \widetilde{\varepsilon},$$

*Then, defining $w_t = \widehat{w}_t / \|\widehat{w}_t\|$ and $w_t^* = \widehat{w}_t^* / \|\widehat{w}_t^*\|$, it satisfies*

$$\|w_t - w_t^*\| \leq \widetilde{\varepsilon} \cdot \Gamma(M, t),$$

*where*

$$\Gamma(M, t) \stackrel{\text{def}}{=} \frac{2}{\lambda_d^t} \begin{cases} t, & \text{if } \lambda_1 = 1; \\ (\lambda_1^t - 1)/(\lambda_1 - 1), & \text{if } \lambda_1 \neq 1. \end{cases} \quad \text{and we have } \Gamma(M, t) \leq 2t \cdot \frac{\max\{1, \lambda_1^t\}}{\lambda_d^t}$$

**Theorem B.3** (Inexact Power Method)**.** *Let $M$ be a PSD matrix with eigenvalues $\lambda_1 \geq \cdots \geq \lambda_d$ and the corresponding eigenvectors $u_1, \ldots, u_d$. With probability at least $1 - p$ it holds that, for every $\varepsilon \in (0, 1)$ and every $t \geq T^{\mathrm{PM}}(\kappa, \varepsilon/4, p)$, if $w_t$ is generated by the power method with per-iteration error $\widetilde{\varepsilon} = \frac{\varepsilon}{4\Gamma(M, t)}$, then*

$$\sum_{i \in [d], \lambda_i \leq (1 - 1/\kappa)\lambda_1} (w_t^\top u_i)^2 \leq \varepsilon \quad and \quad w_t^\top M w_t \geq (1 - 1/\kappa - \varepsilon)\lambda_1 \ .$$

## B.2 Proof of Theorem 3.1

We prove Theorem 3.1 by first showing the following lemma. Most of these properties are analogous to their original variants in [13, 14], but here we take extra care also on negative eigenvalues and thus allowing $M$ to be non-PSD.

**Lemma B.4** (useful properties of `AppxPCA`$^\pm$)**.** *With probability at least $1 - p$, it holds that (by letting $\lambda^* = \|M\|_2$):*

*(a) $\widetilde{\varepsilon}_1 \leq \frac{1}{32\Gamma((\lambda^{(s-1)}I - M)^{-1}, m_1)}$ and $\widetilde{\varepsilon}_1 \leq \frac{1}{32\Gamma((\lambda^{(s-1)}I + M)^{-1}, m_1)}$ for each iteration $s \geq 1$;*

*(b) $\widetilde{\varepsilon}_2 \leq \frac{\varepsilon}{4\Gamma((\lambda^{(f)}I - M)^{-1}, m_2)}$ and $\widetilde{\varepsilon}_2 \leq \frac{\varepsilon}{4\Gamma((\lambda^{(f)}I + M)^{-1}, m_2)}$ when the repeat-until loop is over;*

*(c) $0 \leq \frac{3}{4}(\lambda^{(s-1)} - \lambda^*) \leq \Delta^{(s)} \leq \lambda^{(s-1)} - \lambda^*$ and $\frac{1}{2}(\lambda^{(s-1)} - \lambda^*) \leq \lambda^{(s)} - \lambda^*$ for each iteration $s \geq 1$; and*

*(d) $\lambda^{(f)} - \lambda^* \in [\frac{\delta_\times}{48}\lambda^{(f)}, \frac{\delta_\times}{13}\lambda^*]$ when the repeat-until loop is over.*

*(e) when the repeat-until loop is over,*

$$\text{if } w_a^\top v_a \geq w_b^\top v_b \text{ then } \lambda^{(f)} - \lambda_{\max}(M) \leq \frac{10}{3}(\lambda^{(f)} - \lambda^*); \text{ or}$$
$$\text{if } w_a^\top v_a \leq w_b^\top v_b \text{ then } \lambda^{(f)} + \lambda_{\min}(M) \leq \frac{10}{3}(\lambda^{(f)} - \lambda^*) \ .$$

*Proof.* We denote by $C^{(s)} \stackrel{\text{def}}{=} (\lambda^{(s)}I - M)^{-1}$ and by $D^{(s)} \stackrel{\text{def}}{=} (\lambda^{(s)}I + M)^{-1}$ for notational simplicity. Below we prove all the items by induction for a specific iteration $s \geq 2$ assuming that the items of the previous $s - 1$ iterations are true. The base case of $s = 1$ can be verified similar to the general arguments after some notational changes. We omitted the proofs of the base case $s = 1$.

(a) Recall that

$$\Gamma(C^{(s-1)}, t) \leq 2t \cdot \frac{\max\{1, \lambda_{\max}(C^{(s-1)})^t\}}{\lambda_{\min}(C^{(s-1)})^t} \quad \text{and} \quad \Gamma(D^{(s-1)}, t) \leq 2t \cdot \frac{\max\{1, \lambda_{\max}(D^{(s-1)})^t\}}{\lambda_{\min}(D^{(s-1)})^t}$$

On one hand, we have $\lambda_{\max}(C^{(s-1)}) = \frac{1}{\lambda^{(s-1)} - \lambda^*} \leq \frac{2}{\lambda^{(s-2)} - \lambda^*} \leq \frac{2}{\Delta^{(s-1)}}$ using Lemma B.4.c of the previous iteration. Combining this with the termination criterion $\Delta^{(s-1)} \geq \frac{\delta_\times}{12}\lambda^{(s-1)}$, we have $\lambda_{\max}(C^{(s-1)}) \leq \frac{24}{\delta_\times \lambda^{(s-1)}}$. On the other hand, we have $\lambda_{\min}(C^{(s-1)}) = \frac{1}{\lambda^{(s-1)} - \lambda_{\min}(M)} \geq \frac{1}{\lambda^{(s-1)} + \lambda^*} \geq \frac{1}{2\lambda^{(s-1)}}$. Combining the two bounds we conclude that $\Gamma(C^{(s-1)}, t) \leq 2t(48/\delta_\times)^t$. It is now obvious that $\widetilde{\varepsilon}_1 \leq \frac{1}{32\Gamma(C^{(s-1)}, m_1)}$ is satisfied because $\widetilde{\varepsilon}_1 = \frac{1}{64m_1}(\frac{\delta_\times}{48})^{m_1}$.

Similarly, on one hand, we have $\lambda_{\max}(D^{(s-1)}) = \frac{1}{\lambda^{(s-1)} + \lambda_{\min}(M)} \leq \frac{1}{\lambda^{(s-1)} - \lambda^*} \leq \frac{2}{\lambda^{(s-2)} - \lambda^*} \leq \frac{2}{\Delta^{(s-1)}}$ using Lemma B.4.c of the previous iteration. Combining this with the termination criterion $\Delta^{(s-1)} \geq \frac{\delta_\times}{12}\lambda^{(s-1)}$, we have $\lambda_{\max}(D^{(s-1)}) \leq \frac{24}{\delta_\times \lambda^{(s-1)}}$. On the other hand, we have $\lambda_{\min}(D^{(s-1)}) = \frac{1}{\lambda^{(s-1)} + \lambda_{\max}(M)} \geq \frac{1}{\lambda^{(s-1)} + \lambda^*} \geq \frac{1}{2\lambda^{(s-1)}}$. Combining the two bounds we conclude that $\Gamma(D^{(s-1)}, t) \leq 2t(48/\delta_\times)^t$. It is now obvious that $\widetilde{\varepsilon}_1 \leq \frac{1}{32\Gamma(D^{(s-1)}, m_1)}$ is satisfied.

(b) The same analysis as in the proof of Lemma B.4.a suggests that $\Gamma(C^{(f)}, t) \leq 2t(48/\delta_\times)^t$ and $\Gamma(D^{(f)}, t) \leq 2t(48/\delta_\times)^t$. These immediately imply $\widetilde{\varepsilon}_2 \leq \frac{\varepsilon}{4\Gamma(C^{(f)}, m_2)}$ and $\widetilde{\varepsilon}_2 \leq \frac{\varepsilon}{4\Gamma(D^{(f)}, m_2)}$ because $\widetilde{\varepsilon}_2 = \frac{\varepsilon}{8m_2}(\frac{\delta_\times}{48})^{m_2}$

(c) Because Lemma B.4.a holds for the current iteration $s$ we can apply Theorem B.3 (with $\varepsilon = 1/16$ and $\kappa = 16$) and get

$$w_a^\top C^{(s-1)} w_a \geq \frac{7}{8}\lambda_{\max}(C^{(s-1)}) \quad \text{and} \quad w_b^\top D^{(s-1)} w_b \geq \frac{7}{8}\lambda_{\max}(D^{(s-1)}) \ .$$

By the definition of $v$ in $\texttt{AppxPCA}^\pm$ and the Cauchy-Schwartz inequality it holds that

$$w_a^\top v_a = w_a^\top C^{(s-1)} w_a + w_a^\top (v_a - C^{(s-1)} w_a) \in \left[ w_a^\top C^{(s-1)} w_a - \widetilde{\varepsilon}_1, w_a^\top C^{(s-1)} w_a + \widetilde{\varepsilon}_1 \right] \ , \quad \text{and}$$
$$w_b^\top v_b = w_b^\top D^{(s-1)} w_b + w_b^\top (v_b - D^{(s-1)} w_b) \in \left[ w_b^\top D^{(s-1)} w_b - \widetilde{\varepsilon}_1, w_b^\top D^{(s-1)} w_b + \widetilde{\varepsilon}_1 \right] \ .$$

Combining the above equations we have

$$
\begin{aligned}
w_a^\top v_a - \widetilde{\varepsilon}_1 &\in \left[ \frac{7}{8}\lambda_{\max}(C^{(s-1)}) - 2\widetilde{\varepsilon}_1, \lambda_{\max}(C^{(s-1)}) \right] \\
&\subseteq \left[ \frac{3}{4}\lambda_{\max}(C^{(s-1)}), \lambda_{\max}(C^{(s-1)}) \right] = \left[ \frac{3}{4}, 1 \right] \cdot \frac{1}{\lambda^{(s-1)} - \lambda_{\max}(M)} \ , \quad and \\
w_b^\top v_b - \widetilde{\varepsilon}_1 &\subseteq \left[ \frac{3}{4}\lambda_{\max}(D^{(s-1)}), \lambda_{\max}(D^{(s-1)}) \right] = \left[ \frac{3}{4}, 1 \right] \cdot \frac{1}{\lambda^{(s-1)} + \lambda_{\min}(M)} \ .
\end{aligned}
\tag{B.1}
$$

In other words, $\Delta^{(s)} \stackrel{\text{def}}{=} \frac{3}{4} \cdot \frac{1}{\max\{w_a^\top v_a, w_b^\top w_b\} - \widetilde{\varepsilon}_1} \in \left[ \frac{3}{4}(\lambda^{(s-1)} - \lambda^*), \lambda^{(s-1)} - \lambda^* \right]$ because $\lambda^* = \max\{\lambda_{\max}(M), -\lambda_{\min}(M)\}$.

At the same time, our update rule $\lambda^{(s)} = \lambda^{(s-1)} - \Delta^{(s)}/2$ ensures that $\lambda^{(s)} - \lambda^* = \lambda^{(s-1)} - \lambda^* - \Delta^{(s)}/2 \geq \lambda^{(s-1)} - \lambda^* - \frac{\lambda^{(s-1)} - \lambda^*}{2} = \frac{1}{2}(\lambda^{(s-1)} - \lambda^*)$.

(d) The upper bound holds because $\lambda^{(f)} - \lambda^* = \lambda^{(f-1)} - \frac{\Delta^{(f)}}{2} - \lambda^* \leq \left( \frac{4}{3} - \frac{1}{2} \right)\Delta^{(f)} \leq \frac{5\delta_\times \lambda^{(f)}}{72}$ where the first inequality follows from Lemma B.4.c of this last iteration, and the second inequality follows from our termination criterion $\Delta^{(f)} \leq \frac{\delta_\times \lambda^{(f)}}{12}$. Simply rewriting this inequality we have $\lambda^{(f)} - \lambda^* \leq \frac{5\delta_\times/72}{1 - 5\delta_\times/72}\lambda^* < \frac{\delta_\times}{13}\lambda^*$.

The lower bound is because using Lemma B.4.c (of this and the previous iteration) we have $\lambda^{(f)} - \lambda^* \geq \frac{1}{4}\left( \lambda^{(f-2)} - \lambda^* \right) \geq \frac{\Delta^{(f-1)}}{4} \overset{\text{①}}{\geq} \frac{\delta_\times \lambda^{(f-1)}}{48} \geq \frac{\delta_\times \lambda^{(f)}}{48}$. Here, inequality ① is because $\Delta^{(f-1)} > \frac{\delta_\times \lambda^{(f-1)}}{12}$ due to the termination criterion.

(e) We only prove the case when $w_a^\top v_a \geq w_b^\top v_b$ and the other case is similar. We compute that

$$\lambda^{(f)} - \lambda_{\max}(M) = \lambda^{(f-1)} - \lambda_{\max}(M) - \frac{\Delta^{(f)}}{2} \overset{①}{\leq} \frac{4}{3}(\lambda^{(f-1)} + \lambda_{\min}(M)) - \frac{\Delta^{(f)}}{2}$$

$$= \frac{4}{3}(\lambda^{(f)} + \lambda_{\min}(M)) + \frac{\Delta^{(f)}}{2} \overset{②}{\leq} \frac{4}{3}(\lambda^{(f)} + \lambda_{\min}(M)) + \frac{\delta_\times \lambda^{(f)}}{24}$$

$$\overset{③}{\leq} \frac{4}{3}(\lambda^{(f)} + \lambda_{\min}(M)) + 2(\lambda^{(f)} - \lambda^*) \overset{④}{\leq} \frac{10}{3}(\lambda^{(f)} - \lambda^*) \ .$$

Above, ① is from (B.1) together with the fact that $w_a^\top v_a \geq w_b^\top v_b$; ② is using the termination criterion $\Delta^{(f)} \leq \frac{\delta_\times \lambda^{(f)}}{12}$; ③ is from Lemma B.4.d

Finally, since the success of Theorem B.3 only depends on the randomness of $\widehat{w}_0$, we have that with probability at least $1 - p$ all the above items are satisfied. □

We are now ready to prove Theorem 3.1.

*Proof of Theorem 3.1.* We only focus on the case when $sgn = +$ and the other case is similar. It follows from Theorem B.3 (with $\kappa = 2$) that, letting $\mu_i = 1/(\lambda^{(f)} - \lambda_i)$ be the $i$-th largest eigenvalue of the matrix $(\lambda^{(f)}I - M)^{-1}$, then

$$\sum_{i \in [d], \mu_i \leq \mu_1/2} (w^\top u_i)^2 \leq \varepsilon \ .$$

Note that if an index $i \in [d]$ satisfies $\lambda^* - \lambda_i \geq \frac{\delta_\times}{2}\lambda^*$, then we must have $\lambda^* - \lambda_i \geq \frac{13}{2}(\lambda^{(f)} - \lambda^*)$ owing to $\lambda^{(f)} - \lambda^* \leq \frac{\delta_\times}{13}\lambda^*$ from Lemma B.4.d. This further implies that $\lambda^{(f)} - \lambda_i \geq \frac{15}{2}(\lambda^{(f)} - \lambda^*)$. Plugging in Lemma B.4.e we further have $\lambda^{(f)} - \lambda_i \geq \frac{15}{2} \cdot \frac{3}{10}(\lambda^{(f)} - \lambda_1) > 2(\lambda^{(f)} - \lambda_1)$. Using the definition of $\mu_i$, we must have $\mu_1/2 > \mu_i$. In sum, we also have

$$\sum_{i \in [d], \lambda_i \leq (1 - \delta_\times/2)\lambda^*} (w^\top u_i)^2 \leq \varepsilon \ .$$

On the other hand,

$$w^\top M w = \sum_{i=1}^d \lambda_i (w^\top u_i)^2 \geq -\varepsilon\lambda^* + \sum_{i \in [d], \lambda_i > (1-\delta_\times/2)\lambda^*} \lambda_i (w^\top u_i)^2$$

$$\geq -\varepsilon\lambda^* + (1 - \delta_\times/2)\lambda^* \cdot \sum_{i \in [d], \lambda_i > (1-\delta_\times/2)\lambda^*} (w^\top u_i)^2$$

$$\geq -\varepsilon\lambda^* + (1 - \delta_\times/2)(1 - \varepsilon)\lambda^* \geq (1 - \delta_\times/2)(1 - 3\varepsilon)\lambda^* \ .$$

The number of oracle calls to $\mathcal{A}$ is determined by the number of iterations in the repeat-until loop. It is easy to verify that there are at most $O(\log(1/\delta_\times))$ such iteartions, so the total number of oracle calls to $\mathcal{A}$ is only $O(\log(1/\delta_\times)m_1 + m_2)$.

As for the condition number, each time we call $\mathcal{A}$ we have

$$\frac{\lambda_{\max}(\lambda^{(s)}I - M)}{\lambda_{\min}(\lambda^{(s)}I - M)} \leq \frac{\lambda^{(s)} - \lambda_d}{\lambda^{(s)} - \lambda_1} \leq \frac{\lambda^{(s)} + \lambda^*}{\lambda^{(s)} - \lambda^*} \leq \frac{2\lambda^{(s)}}{\lambda^{(s)} - \lambda^*} \quad \text{and} \quad \frac{\lambda_{\max}(\lambda^{(s)}I + M)}{\lambda_{\min}(\lambda^{(s)}I + M)} \leq \frac{2\lambda^{(s)}}{\lambda^{(s)} - \lambda^*}$$

If $s = 0$ then we have $\frac{\lambda^{(0)}}{\lambda^{(0)} - \lambda^*} \leq \frac{1 + \delta_\times}{\delta_\times}$ because $\lambda^* \leq 1$. If $s \leq f - 2$ then we have $\frac{\lambda^{(s)}}{\lambda^{(s)} - \lambda^*} \leq \frac{\lambda^{(s)}}{\Delta^{(s+1)}} \leq \frac{\lambda^{(s)}}{\delta_\times \lambda^{(s+1)}/12} \leq \frac{12}{\delta_\times}$ where the first inequality follows from Lemma B.4.c, the second inequality follows from the stopping criterion, and the third inequality follows from the monotonicity of $\lambda^{(s)}$. If $s = f - 1$ then we have $\frac{\lambda^{(s)}}{\lambda^{(s)} - \lambda^*} \leq \frac{2\lambda^{(s)}}{\lambda^{(s-1)} - \lambda^*} \leq \frac{2\lambda^{(s)}}{\Delta^{(s)}} \leq \frac{2\lambda^{(s)}}{\delta_\times \lambda^{(s)}/12} = \frac{24}{\delta_\times}$ where the first two inequalities follow from Lemma B.4.c and the third inequality follows from our stopping criterion. If $s = f$ then we have $\frac{\lambda^{(s)}}{\lambda^{(s)} - \lambda^*} \leq \frac{48}{\delta_\times}$ owing to Lemma B.4.d. In all cases we have $\frac{\lambda_{\max}(\lambda^{(s)}I - M)}{\lambda_{\min}(\lambda^{(s)}I - M)} \leq \frac{96}{\delta_\times}$ and $\frac{\lambda_{\max}(\lambda^{(s)}I + M)}{\lambda_{\min}(\lambda^{(s)}I + M)} \leq \frac{96}{\delta_\times}$.

Finally, we have $\frac{1}{\lambda_{\min}(\lambda^{(s)}I - M)} = \frac{\lambda^{(s)}}{\lambda_{\min}(\lambda^{(s)}I - M)} \cdot \frac{1}{\lambda^{(s)}} \leq \frac{\lambda^{(s)}}{\lambda^{(s)} - \lambda^*} \cdot \frac{1}{\lambda^*} \leq \frac{48}{\delta_\times \lambda^*}$ where the first inequality follows from $\lambda^{(s)} \geq \lambda^*$. Similarly, we also have $\frac{1}{\lambda_{\min}(\lambda^{(s)}I + M)} \leq \frac{48}{\delta_\times \lambda^*}$. $\qquad\square$

# C   Main Matrix-Algebra Lemmas

In this section we provide some necessary lemmas on matrix algebra that shall become essential for our proof of Theorem 4.1. Many of these lemmas are analogous to those ones used in the SVD algorithm by the same authors of this paper [2], however, we need some extra care in this paper because the underlying matrix $M$ is no longer PSD.

**Proposition C.1.** *Let $A, B$ be two (column) orthonormal matrix such that for $\eta \geq 0$,*
$$A^\top B B^\top A \succeq (1 - \eta)I$$
*Then we have: there exists a matrix $Q, \|Q\|_2 \leq 1$ such that*
$$\|A - BQ\|_2 \leq \sqrt{\eta}$$

*Proof.* Since $A^\top A = I$ and $A^\top B B^\top A \succeq (1 - \eta)I$, we know that $A^\top B^\perp (B^\perp)^\top A \preceq \eta I$. By the fact that
$$A = (BB^\top + B^\perp (B^\perp)^\top)A = BB^\top A + B^\perp (B^\perp)^\top A$$
we can let $Q = B^\top A$ and obtain
$$\|A - BQ\|_2 \leq \|B^\perp (B^\perp)^\top A\|_2 \leq \sqrt{\eta} \ . \qquad\square$$

## C.1   Approximate Projection Lemma

The next lemma states that, projecting a symmetric matrix $M$ into the orthogonal space of $V_s \in \mathbb{R}^{d \times s}$ is almost equivalent to projecting it into the orthogonal space of $Q_s \in \mathbb{R}^{d \times s}$, if $Q_s$ is the projection of $V_s$ into the orthogonal space of $U$ but $\|V_s^\top U\|$ is small. This lemma is obvious if "small" means zero correlation: if $V_s$ were completely orthogonal to $U$ then $Q_s$ would equal to $V_s$, so projecting $M$ into the orthogonal space of $V_s$ would be equivalent to that of $Q_s$. However, even in the inexact scenario, this argument is true.

**Lemma C.2.** *Let $M$ be a symmetric matrix with (not necessarily sorted) eigenvalues $\lambda_1, \ldots, \lambda_d$ and the corresponding (normalized) eigenvectors $u_1, \ldots, u_d \in \mathbb{R}^d$. For every $k \geq 1$, define $U^\perp = (u_1, \ldots, u_k) \in \mathbb{R}^{d \times k}$ and $U = (u_{k+1}, \ldots, u_d) \in \mathbb{R}^{d \times (d-k)}$. For every $\varepsilon \in (0, \frac{1}{2})$, let $V_s \in \mathbb{R}^{d \times s}$ be a column orthogonal matrix such that $\|V_s^\top U\|_2 \leq \varepsilon$, define $Q_s \in \mathbb{R}^{d \times s}$ to be an arbitrary orthogonal basis of the column span of $U^\perp (U^\perp)^\top V_s$, then we have:*
$$\left\| \left(I - Q_s Q_s^\top\right) M \left(I - Q_s Q_s^\top\right) - \left(I - V_s V_s^\top\right) M \left(I - V_s V_s^\top\right) \right\|_2 \leq 13\varepsilon \|M\|_2 \ .$$

*Proof of Lemma C.2.* Since $Q_s$ is an orthogonal basis of the column span of $U^\perp (U^\perp)^\top V_s$, there is a matrix $R \in \mathbb{R}^{s \times s}$ such that
$$Q_s = U^\perp (U^\perp)^\top V_s R$$

Using the fact that $Q_s^\top Q_s = I$, we have:
$$(U^\perp (U^\perp)^\top V_s R)^\top (U^\perp (U^\perp)^\top V_s R) = I \implies R^\top V_s^\top U^\perp (U^\perp)^\top V_s R = I \ .$$

By the fact that $V_s^\top V_s = I$ and $U^\perp (U^\perp)^\top + UU^\top = I$, we can rewrite the above equality as:

$$R^\top \left( I - V_s^\top UU^\top V_s \right) R = I \tag{C.1}$$

From our lemma assumption, we have: $\|V_s^\top U\|_2 \le \varepsilon$, which implies $0 \preceq V_s^\top UU^\top V_s \preceq \varepsilon^2 I$. Putting this into (C.1), we obtain:

$$I \preceq R^\top R \preceq \frac{1}{1-\varepsilon^2} I \preceq \left(1 + \frac{4}{3}\varepsilon^2\right) I$$

The above inequality directly implies that $I \preceq RR^\top \preceq \left(1 + \frac{4}{3}\varepsilon^2\right) I$. Therefore,

$$
\begin{aligned}
& \left\| Q_s Q_s^\top - V_s V_s^\top \right\|_2 \\
= \ & \left\| U^\perp (U^\perp)^\top V_s RR^\top V_s^\top U^\perp (U^\perp)^\top - V_s V_s^\top \right\|_2 \\
= \ & \left\| U^\perp (U^\perp)^\top V_s RR^\top V_s^\top U^\perp (U^\perp)^\top - (U^\perp (U^\perp)^\top + UU^\top) V_s V_s^\top (U^\perp (U^\perp)^\top + UU^\top) \right\|_2 \\
\le \ & \left\| U^\perp (U^\perp)^\top V_s (RR^\top - I) V_s^\top U^\perp (U^\perp)^\top \right\|_2 + \left\| UU^\top V_s V_s^\top UU^\top \right\|_2 + 2 \left\| U^\perp (U^\perp)^\top V_s V_s^\top UU^\top \right\|_2 \\
\le \ & \left\| RR^\top - I \right\|_2 + \left\| U^\top V_s V_s^\top U \right\|_2 + 2 \left\| V_s^\top UU^\top V_s \right\|_2^{1/2} \\
\le \ & \frac{4}{3}\varepsilon^2 + \varepsilon^2 + 2\varepsilon < \frac{19}{6}\varepsilon \ .
\end{aligned}
$$

Finally, we have

$$
\begin{aligned}
& \left\| \left( I - Q_s Q_s^\top \right) M \left( I - Q_s Q_s^\top \right) - \left( I - V_s V_s^\top \right) M \left( I - V_s V_s^\top \right) \right\|_2 \\
\le \ & 2 \left\| \left( Q_s Q_s^\top - V_s V_s^\top \right) M \right\|_2 + \left\| \left( Q_s Q_s^\top - V_s V_s^\top \right) M Q_s Q_s^\top \right\|_2 + \left\| \left( Q_s Q_s^\top - V_s V_s^\top \right) M V_s V_s^\top \right\|_2 \\
\le \ & \frac{19 \times 4}{6}\varepsilon \|M\|_2 < 13\varepsilon \|M\|_2 \ . \qquad \square
\end{aligned}
$$

## C.2 Gap-Free Wedin Theorem

**Lemma C.3** (two-sided gap-free Wedin theorem). *For $\varepsilon \ge 0$, let $A, B$ be two symmetric matrices such that $\|A - B\|_2 \le \varepsilon$. For every $\mu \ge 0$, $\tau > 0$, let $U$ be column orthonormal matrix consisting of eigenvectors of $A$ with **absolute** eigenvalues $\le \mu$, let $V$ be column orthonormal matrix consisting of eigenvectors of $B$ with **absolute** eigenvalues $\ge \mu + \tau$, then we have:*

$$\|U^\top V\| \le \frac{\varepsilon}{\tau} \ .$$

*Proof of Lemma C.3.* We write $A$ and $B$ in terms of eigenvalue decomposition:

$$A = U\Sigma U^\top + U'\Sigma' U'^\top \quad \text{and} \quad B = V\widetilde{\Sigma} V^\top + V'\widetilde{\Sigma}' V'^\top \ ,$$

where $U'$ is orthogonal to $U$ and $V'$ is orthogonal to $V$. Letting $R = A - B$, we obtain:

$$\Sigma U^\top = U^\top A = U^\top (B + R)$$
$$\Longrightarrow \Sigma U^\top V = U^\top BV + U^\top RV = U^\top V\widetilde{\Sigma} + U^\top RV$$
$$\Longrightarrow \Sigma U^\top V \widetilde{\Sigma}^{-1} = U^\top V + U^\top RV \widetilde{\Sigma}^{-1} \ .$$

Taking spectral norm on both sides, we obtain:

$$\|\Sigma\|_2 \|U^\top V\|_2 \|\widetilde{\Sigma}^{-1}\|_2 \ge \|\Sigma U^\top V \widetilde{\Sigma}^{-1}\|_2 \ge \|U^\top V\|_2 - \|U^\top RV \widetilde{\Sigma}^{-1}\|_2 \ .$$

This can be simplified to

$$\frac{\mu}{\mu + \tau} \|U^\top V\|_2 \ge \|U^\top V\|_2 - \frac{\varepsilon}{\mu + \tau} \ ,$$

and therefore we have $\|U^\top V\|_2 \leq \frac{\varepsilon}{\tau}$ as desired. $\qquad\square$

### C.3 Eigenvector Projection Lemma

Our next technical lemma studies the projection of a matrix $M$ into the orthogonal direction of a vector $v$, where $v$ has little correlation with $M$'s leading eigenvectors below some threshold $\mu$ (denoted by $U$). The conclusion of the lemma says that, after the projection, if we study the leading eigenvectors of $M' = (I - vv^\top)M(I - vv^\top)$ below some threshold $\mu + \tau$ and denote it by $V_1$, then $U$ approximately embeds into $V_1$, meaning that although $V_1$ could be of a larger dimension of $U$, however, there exists a matrix $Q$ with spectral norm no more than 1 such that $\|U - V_1 Q\|$ is small.

**Lemma C.4.** *Let $M \in \mathbb{R}^{d\times d}$ be a symmetric matrix with eigenvalues $\lambda_1, \dots, \lambda_d$ and corresponding eigenvectors $u_1, \dots, u_d$. Suppose $|\lambda_1| \geq \cdots \geq |\lambda_d|$. Define $U = (u_{j+1}, \dots, u_d) \in \mathbb{R}^{d\times(d-j)}$ to be the matrix consisting of all eigenvectors with* **absolute** *eigenvalues $\leq \mu$. Let $v \in \mathbb{R}^d$ be a unit vector such that $\|v^\top U\|_2 \leq \varepsilon \leq 1/2$, and define*

$$M' = \left(I - vv^\top\right) M \left(I - vv^\top\right) \ .$$

*Then, denoting by $[V_2, V_1, v] \in \mathbb{R}^{d\times d}$ the unitary matrix consisting of (column) eigenvectors of $M'$, where $V_1$ consists of eigenvectors with absolute eigenvalue $\leq \mu + \tau$, then there exists a matrix $Q$ with spectral norm $\|Q\|_2 \leq 1$ such that*

$$\|U - V_1 Q\|_2 \leq \sqrt{\frac{169\varepsilon^2\|M\|_2^2}{\tau^2} + \varepsilon^2} \ .$$

*Proof of Lemma C.4.* Using Lemma C.2, let $q = \frac{U^\perp(U^\perp)^\top v}{\|U^\perp(U^\perp)^\top v\|_2}$ be the projection of $v$ to $U^\perp$, we know that

$$\left\| \left(I - qq^\top\right) M \left(I - qq^\top\right) - \left(I - vv^\top\right) M \left(I - vv^\top\right) \right\|_2 \leq 13\varepsilon\|M\|_2 \ .$$

Denote $\left(I - qq^\top\right) M \left(I - qq^\top\right)$ as $M''$. We know that $u_{j+1}, \dots, u_d$ are still eigenvectors of $M''$ with eigenvalue $\lambda_{j+1}, \dots, \lambda_d$.

Apply Lemma C.3 on $A = M''$, $U$ and $B = M'$, $V = V_2$, we obtain:

$$\|U^\top V_2\|_2 \leq \frac{13\varepsilon\|M\|_2}{\tau} \ .$$

This implies that

$$U^\top V_1 V_1^\top U = I - U^\top V_2 V_2^\top U - U^\top vv^\top U \succeq \left(1 - \frac{169\varepsilon^2\|M\|_2^2}{\tau^2} - \varepsilon^2\right) I \ ,$$

where the inequality uses the assumption $\|v^\top U\|_2 \leq \varepsilon$.

Apply Proposition C.1 to $A = U$ and $B = V_1$, we conclude that there exists a matrix $Q$, $\|Q\|_2 \leq 1$ such that

$$\|U - V_1 Q\|_2 \leq \sqrt{\frac{169\varepsilon^2\|M\|_2^2}{\tau^2} + \varepsilon^2} \ . \qquad\square$$

## D  Matrix Inversion via Approx Accelerated Gradient Descent

Given a positive definite matrix $\mathbb{N}$, it is well-known that one can reduce the (approximate) matrix inversion problem $\mathbb{N}^{-1}\chi$ to multiple computations of the matrix-vector multiplication (i.e., of the form $w' \leftarrow \mathbb{N}w$). In particular, Chebyshev method [9] uses the so-called Chebyshev polynomial for

---

**Algorithm 4** $\mathtt{AGD^{inexact}}(f, x_0, T)$

---

**Input:** $f$ an $L$-smooth and $\sigma$-strongly convex function;
   $x_0$ some initial point; and
   $T$ the number of iterations.

**Output:** $y_T$.

1:  $\tau \leftarrow \frac{2}{1+\sqrt{8L/\sigma+1}}$, $\eta \leftarrow \frac{1}{\tau L}$.         $\diamond$   $\tau = O(\frac{\sqrt{\sigma}}{\sqrt{L}})$ and $\eta = O(\frac{1}{\sqrt{\sigma L}})$

2:  $y_0 \leftarrow x_0,$  $z_0 \leftarrow x_0$.

3:  **for** $k \leftarrow 0$ **to** $T-1$ **do**

4:   $x_{k+1} \leftarrow \tau z_k + (1-\tau)y_k$.

5:   Compute approximate gradient $\widetilde{\nabla}f(x_{k+1})$ satisfying $\|\widetilde{\nabla}f(x_{k+1}) - \nabla f(x_{k+1})\|_2 \leq \widetilde{\varepsilon}$.

6:   $y_{k+1} \leftarrow x_{k+1} - \frac{1}{L}\widetilde{\nabla}f(x_{k+1})$

7:   $z_{k+1} \leftarrow \frac{1}{1+\eta\sigma}(z_k + \eta\sigma x_{k+1} - \eta\widetilde{\nabla}f(x_{k+1}))$

8:  **end for**

9:  **return** $y_T$.

---

this purpose, and the number of matrix-vector multiplications is determined by the degree of that polynomial.

In this section, we revisit this problem by allowing matrix-vector multiplications to be computed only approximately. We emphasize that this is not a simple task in general. If matrix inversion is reduced to $T$ matrix-vector multiplications, then a standard analysis implies that each of these multiplications must be computed up to a very small error $2^{-\Omega(T)}$. If the actual matrix-vector multiplication subroutine has a logarithmic dependency on the error in its running time, then we will have a total running time at least quadratically dependent on $T$.[9]

To avoid such an exponentially accuracy loss, we abandon known results (such as Chebyshev method) and design our own method. We prove the following theorem in this section:

---

**Theorem D.1.** *Given a matrix* $\mathbb{N}$*, we reduce the problem of computing* $\xi \leftarrow \mathbb{N}^{-1}\chi$ *to multiple computations* $w' \leftarrow \mathbb{N}w$ *as follows.*

*If* $\mathbb{N}$ *satisfies* $\mathbb{N} = B^{-1/2}NB^{1/2}$ *where* $N$ *and* $B$ *are both* $d \times d$ *positive definite matrices, for every* $\widetilde{\varepsilon} > 0$ *and* $\chi \in \mathbb{R}^d$*, in order to obtain* $\xi$ *satisfying* $\|\xi - \mathbb{N}^{-1}\chi\| \leq \widetilde{\varepsilon}$*,*

- *it suffices to compute* $w' \leftarrow \mathbb{N}w$ *only* $\widetilde{O}(\sqrt{\kappa_N})$ *times, and*

- *each time of accuracy* $\|w' - \mathbb{N}w\| \leq O(1/\mathsf{poly}(\kappa_B, \widetilde{\varepsilon}, \lambda_{\min}(N)))$*.*

---

Our reduction is based on an inexact variant of the accelerated gradient descent (AGD) method originally put forward by Nesterov [23], which relies on some convex optimization techniques and can be proved using the linear-coupling framework [5]. We prove this inexact AGD result in Appendix D.1. Our final proof of Theorem D.1 is included in Appendix D.2.

## D.1 Inexact Accelerated Gradient Descent

We study an inexact version of the classical accelerated gradient descent (AGD) method, and our pseudocode is presented in Algorithm 4. The difference between our method and known AGD methods is that we only require the algorithm to know an approximate gradient $\widetilde{\nabla}f(x_{k+1})$ in each iteration $k$, as opposed to the exact full gradient $\nabla f(x_{k+1})$. We require $\|\widetilde{\nabla}f(x_{k+1}) - \nabla f(x_{k+1})\|_2$ to be upper bounded by some parameter $\widetilde{\varepsilon}$ in each iteration. Our next convergence theorem states

---

[9]Indeed, for instance in the ALS algorithm of [29] for solving CCA, the authors obtained a running time proportional to $1/\mathtt{gap}^2$ although there are only $1/\mathtt{gap}$ iterations.

that this inexact AGD method only incurs an additive loss proportional to $O(\widetilde{\varepsilon}^2)$.

**Theorem D.2** (inexact AGD). *If $f(x)$ is $L$-smooth and $\sigma$-strongly convex, then* $\mathtt{AGD}^{\mathsf{inexact}}(f, x_0, T)$ *produces an output $y_T$ satisfying*

$$f(y_T) - f(x^*) \leq O(1) \cdot (1-\tau)^T (f(x_0) - f(x^*)) + O(\frac{\widetilde{\varepsilon}^2}{\sigma}) \ ,$$

*where $\tau = \Omega(\sqrt{\sigma/L})$. In other words, if the approximate gradient oracle satisfies $\widetilde{\varepsilon} \leq O(\sqrt{\varepsilon\sigma})$ and $T = O(\sqrt{L/\sigma} \cdot \log(1/\varepsilon))$, then we have $f(y_T) - f(x^*) \leq \varepsilon$.*

Theorem D.2 can be proved using the linear-coupling framework of [5]. In this framework, accelerated methods are analyzed by a gradient descent lemma (Lemma D.3 below), a mirror descent lemma (Lemma D.4 below), and a coupling step (Lemma D.5 and D.6 below).

**Lemma D.3** (gradient descent). $f(y_{k+1}) \leq f(x_{k+1}) - \frac{1}{2L} \|\nabla f(x_{k+1})\|_2^2 + \frac{\widetilde{\varepsilon}^2}{2L}$.

*Proof.* Abbreviating $x_{k+1}$ by $x$ and $y_{k+1}$ by $y$, the smoothness property of function $f(\cdot)$ tells us

$$f(y) - f(x) \leq \langle \nabla f(x), y - x \rangle + \frac{L}{2} \|y - x\|^2 \ .$$

Now, since $y - x = -\frac{\nabla f(x) + \chi}{L}$ where $\|\chi\|_2 \leq \widetilde{\varepsilon}$, we have

$$\langle \nabla f(x), y - x \rangle + \frac{L}{2} \|y - x\|_2^2 = \frac{-1}{L} \langle \nabla f(x), \nabla f(x) + \chi \rangle + \frac{1}{2L} \langle \nabla f(x) + \chi, \nabla f(x) + \chi \rangle$$

$$\leq -\frac{1}{2L} \|\nabla f(x)\|_2^2 + \frac{\widetilde{\varepsilon}^2}{2L} \ . \qquad \square$$

Since our update on $z$ can be written in the following minimization form, known as mirror-descent form in optimization literatures:

$$z_{k+1}^{(i)} = \min_z \left\{ \frac{1}{2} \|z - z_k\|^2 + \eta \langle \widetilde{\nabla} f(x_{k+1}), z \rangle + \frac{\eta\sigma}{2} \|z - x_{k+1}\|_2^2 \right\} \ . \tag{D.1}$$

It implies the following classical lemma (see for instance [6, Lemma 5.4]):

**Lemma D.4** (mirror descent). *For every $u \in \mathbb{R}^n$,*

$$\eta \langle \widetilde{\nabla} f(x_{k+1}), z_{k+1} - u \rangle - \frac{\eta\sigma}{2} \|x_{k+1} - u\|_2^2 \leq -\frac{1}{2} \|z_k - z_{k+1}\|_2^2 + \frac{1}{2} \|z_k - u\|_2^2 - \frac{1 + \eta\sigma}{2} \|z_{k+1} - u\|_2^2 \ .$$

The following inequality is a nature linear combination of the two lemmas above:

**Lemma D.5** (coupling 1). *For every $u \in \mathbb{R}^n$,*

$$\eta \langle \nabla f(x_{k+1}), z_k - u \rangle - \frac{\eta\sigma}{2} \|u - x_{k+1}\|_2^2$$

$$\leq \eta^2 L \big( f(x_{k+1}) - f(y_{k+1}) \big) + \frac{1}{2} \|z_k - u\|_2^2 - \frac{1 + \eta\sigma/2}{2} \|z_{k+1} - u\|_2^2 + \widetilde{\varepsilon}^2 (\frac{\eta}{\sigma} + \frac{\eta^2}{2}) \ .$$

*Proof.* Combining Lemma D.3 and Lemma D.4 we deduce that for each $i \in [n]$,

$$\langle \eta \nabla f(x_{k+1}), z_k - u \rangle - \frac{\eta\sigma}{2} \|x_{k+1} - u\|_2^2$$

$$\leq \langle \eta \nabla f(x_{k+1}), z_k - z_{k+1} \rangle + \langle \eta \widetilde{\nabla} f(x_{k+1}), z_{k+1} - u \rangle + \widetilde{\varepsilon}\eta \|z_{k+1} - u\|_2 - \frac{\eta\sigma}{2} \|x_{k+1} - u\|_2^2$$

$$\overset{①}{\leq} \langle \eta \nabla f(x_{k+1}), z_k - z_{k+1} \rangle - \frac{1}{2} \|z_k - z_{k+1}\|_2^2 + \frac{1}{2} \|z_k - u\|_2^2 - \frac{1 + \eta\sigma}{2} \|z_{k+1} - u\|_2^2 + \widetilde{\varepsilon}\eta \|z_{k+1} - u\|_2$$

$$\overset{②}{\leq} \frac{\eta^2}{2} \|\nabla f(x_{k+1})\|_2^2 + \frac{1}{2} \|z_k - u\|_2^2 - \frac{1 + \eta\sigma/2}{2} \|z_{k+1} - u\|_2^2 + \frac{\widetilde{\varepsilon}^2\eta}{\sigma}$$

$$\overset{③}{\leq} \eta^2 L\big(f(x_{k+1}) - f(y_{k+1})\big) + \frac{1}{2}\|z_k - u\|_2^2 - \frac{1 + \eta\sigma/2}{2}\|z_{k+1} - u\|_2^2 + \widetilde{\varepsilon}^2(\frac{\eta}{\sigma} + \frac{\eta^2}{2}) \ .$$

Above, ① uses Lemma D.4, ② uses the Young's inequality which states $2\langle a, b\rangle \leq \|a\|^2 + \|b\|^2$, ③ uses Lemma D.3. $\qquad\square$

Taking into account $x_{k+1} = \tau z_k + (1 - \tau)y_k$ and the convexity of $f(\cdot)$, we can rewrite some terms of Lemma D.5 and obtain

**Lemma D.6** (coupling 2)**.**

$$0 \leq \frac{(1-\tau)\eta}{\tau}(f(y_k) - f(x^*)) - \frac{\eta}{\tau}(f(y_{k+1}) - f(x^*)) + \frac{1}{2}\|z_k - x^*\|_2^2 - \frac{1 + \eta\sigma/2}{2}\|z_{k+1} - x^*\|_2^2 + \widetilde{\varepsilon}^2(\frac{\eta}{\sigma} + \frac{\eta^2}{2})$$

*Proof.*

$\eta(f(x_{k+1}) - f(x^*))$

$\overset{①}{\leq} \eta\langle\nabla f(x_{k+1}), x_{k+1} - x^*\rangle - \frac{\eta\sigma}{2}\|x^* - x_{k+1}\|_2^2$

$= \eta\langle\nabla f(x_{k+1}), x_{k+1} - z_k\rangle + \eta\langle\nabla f(x_{k+1}), z_k - x^*\rangle - \frac{\eta\sigma}{2}\|x^* - x_{k+1}\|_2^2$

$\overset{②}{=} \frac{(1-\tau)\eta}{\tau}\langle\nabla f(x_{k+1}), y_k - x_{k+1}\rangle + \eta\langle\nabla f(x_{k+1}), z_k - x^*\rangle - \frac{\eta\sigma}{2}\|x^* - x_{k+1}\|_2^2$

$\overset{③}{\leq} \frac{(1-\tau)\eta}{\tau}(f(y_k) - f(x_{k+1})) + \eta^2 L\big(f(x_{k+1}) - f(y_{k+1})\big) + \frac{1}{2}\|z_k - u\|_2^2 - \frac{1 + \eta\sigma/2}{2}\|z_{k+1} - u\|_2^2 + \widetilde{\varepsilon}^2(\frac{\eta}{\sigma} + \frac{\eta^2}{2}) \ .$

Above, ① is owing to the strong convexity of $f(\cdot)$, ② uses the fact that $x_{k+1} = \tau z_k + (1-\tau)y_k$, and ③ uses the convexity of $f(\cdot)$ as well as Lemma D.5 with the choice of $u = x^*$. Recall $\eta = \frac{1}{\tau L}$, we arrive at the desired inequality. $\qquad\square$

We are now ready to prove Theorem D.2.

*Proof of Theorem D.2.* We choose $\tau = \frac{2}{1 + \sqrt{8L/\sigma + 1}} \in [0, 1)$, and this choice ensures that $1 + \eta\sigma/2 = \frac{1}{1-\tau}$. Under these parameter choices, Lemma D.6 becomes

$$\big(f(y_{k+1}) - f(x^*)\big) + \frac{\tau}{2\eta(1-\tau)}\|z_{k+1} - x^*\|_2^2 \leq (1-\tau)\Big((f(y_k) - f(x^*)) + \frac{\tau}{2\eta(1-\tau)}\|z_k - x^*\|_2^2\Big) + \widetilde{\varepsilon}^2\tau(\frac{1}{\sigma} + \frac{\eta}{2})$$

Telescoping it for all iterations $k = 0, 1, \ldots, T - 1$, we conclude that

$$f(y_T) - f(x^*) \leq (1-\tau)^T\Big(f(y_0) - f(x^*) + \frac{\tau}{2\eta}\|z_0 - x^*\|_2^2\Big) + \widetilde{\varepsilon}^2(\frac{1}{\sigma} + \frac{\eta}{2}) \leq O(1)\cdot(1-\tau)^T(f(x_0) - f(x^*)) + O\big(\frac{\widetilde{\varepsilon}^2}{\sigma}\big) \ .$$

where the last inequality is because (i) $x_0 = y_0 = z_0$, (ii) $\tau/\eta = O(\sigma)$ and (iii) the strong convexity of $f(\cdot)$ which implies $f(x_0) - f(x^*) \geq \frac{\sigma}{2}\|x_0 - x^*\|_2^2$. $\qquad\square$

## D.2 Proof of Theorem D.1

*Proof of Theorem D.1.* We first verify accuracy. Since

$$\|\xi - N^{-1}\chi\| \leq \widetilde{\varepsilon} \Longleftarrow \|B^{1/2}\xi - B^{1/2}N^{-1}\chi\| \leq \widetilde{\varepsilon} \cdot \sqrt{\lambda_{\min}(B)}$$
$$\Longleftarrow \|B^{1/2}\xi - N^{-1}B^{1/2}\chi\| \leq \widetilde{\varepsilon} \cdot \sqrt{\lambda_{\min}(B)} \ , \qquad (D.2)$$

it suffices to find $\xi$ to satisfy (D.2) in order to satisfy the accuracy requirement $\|\xi - N^{-1}\chi\| \leq \widetilde{\varepsilon}$. Define $f(x) \overset{\text{def}}{=} \frac{1}{2}x^\top Nx - \big(B^{1/2}\chi\big)^\top x$ and let $x^*$ be its minimizer. Then it satisfies $x^* = N^{-1}B^{1/2}\chi$

and $f(x) - f(x^*) = \frac{1}{2}(x - x^*)^\top N(x - x^*)$. For this reason, it suffices to find an approximate minimizer of $f(x)$ satisfying

$$\frac{1}{2}(x - x^*)^\top N(x - x^*) = f(x) - f(x^*) \leq \frac{\widetilde{\varepsilon}^2}{2}\lambda_{\min}(B)\lambda_{\min}(N) =: \widetilde{\varepsilon}' \tag{D.3}$$

because if we let $\xi = B^{-1/2}x$ then the above inequality implies $\frac{1}{2}\|x - x^*\|^2 \leq \frac{\widetilde{\varepsilon}^2}{2} \cdot \lambda_{\min}(B)$ which is the same as (D.2). In sum, we can call $\texttt{AGD}^{\texttt{inexact}}$ to find an approximate minimizer $x$ with additive error no more than $\widetilde{\varepsilon}'$, and then defining $\xi = B^{-1/2}x$ gives a solution of $\xi$ satisfying $\|\xi - N^{-1}\chi\| \leq \widetilde{\varepsilon}$.

We now focus on the actual implementation of $\texttt{AGD}^{\texttt{inexact}}$. If we choose $x_0 = 0$ as the initial vector, we can write $x_k, y_k, z_k$ implicitly as $x_k = B^{1/2}\mathbb{x}_k, y_k = B^{1/2}\mathbb{y}_k, z_k = B^{1/2}\mathbb{y}_k$ (thus only keep track of $\mathbb{x}_k, \mathbb{y}_k, \mathbb{z}_k$) throughout the algorithm. Under these notations, we claim that it suffices to perform matrix vector multiplication on $\mathbb{N}$ (i.e., of the form $w' \leftarrow \mathbb{N}w$) for at most $O(T)$ times on those implicit vectors where $T = O\big(\sqrt{\lambda_{\max}(N)/\lambda_{\min}(N)}\log(1/\widetilde{\varepsilon}')\big)$ is the number of iterations of $\texttt{AGD}^{\texttt{inexact}}$ according to Theorem D.2.

This is so because $\nabla f(x_k) = Nx_k - B^{1/2}\chi = B^{1/2}\big(\mathbb{N}\mathbb{x}_k - \chi\big)$ and therefore for instance $y_{k+1} \leftarrow x_{k+1} - \frac{1}{L}\widetilde{\nabla}f(x_{k+1})$ can be implemented as $\mathbb{y}_{k+1} \leftarrow \mathbb{x}_{k+1} - \frac{1}{L}(\mathbb{N}\mathbb{x}_{k+1} - \chi)$ so only matrix-vector multiplication on $\mathbb{N}$ is needed. In addition, as long as each $w' \leftarrow \mathbb{N}w$ is computed to an additive error $\|w' - \mathbb{N}w\| \leq O\big(\widetilde{\varepsilon} \cdot \lambda_{\min}(N)\sqrt{\lambda_{\min}(B)/\lambda_{\max}(B)}\big)$, we can use $B^{1/2}(w' - \chi)$ as the approximate gradient which is different from the true gradient $\nabla f(x_k)$ by an additive amount $O(\sqrt{\lambda_{\min}(N)\widetilde{\varepsilon}'})$. This satisfies the approximation require of Theorem D.2, and thus the accuracy guarantee provided by Theorem D.2 is satisfied. $\qquad\square$

# E  Proof for Section 4: GenEV Theorems

## E.1  Proof of Theorem 4.1

In this section we prove Theorem 4.1 formally.

---

**Theorem 4.1** (restated). *Let $M \in \mathbb{R}^{d \times d}$ be a symmetric matrix with eigenvalues $\lambda_1, \ldots, \lambda_d \in [-1, 1]$ and corresponding eigenvectors $u_1, \ldots, u_d$. Suppose without loss of generality that $|\lambda_1| \geq \cdots \geq |\lambda_d|$.*

*Suppose $k \in [d]$, $\delta_\times, p \in (0,1)$. Then, $\texttt{LazyEV}$ outputs a (column) orthonormal matrix $V_k = (v_1, \ldots, v_k) \in \mathbb{R}^{d \times k}$ which, with probability at least $1 - p$, satisfies all of the following properties. (Denote by $M_k = (I - V_kV_k^\top)M(I - V_kV_k^\top)$.)*

(a) *Core lemma: if $\varepsilon_{\mathsf{pca}} \leq \frac{\varepsilon^4 \delta_\times}{2^{12}k^3(|\lambda_1|/|\lambda_k|)^2}$, then $\|V_k^\top U\|_2 \leq \varepsilon$, where $U = (u_j, \ldots, u_d)$ is the (column) orthonormal matrix and $j$ is the smallest index satisfying $|\lambda_j| \leq (1 - \delta_\times)\|M_{k-1}\|_2$.*

(b) *Spectral norm guarantee: if $\varepsilon_{\mathsf{pca}} \leq \frac{\delta_\times^5}{2^{28}k^3(|\lambda_1|/|\lambda_{k+1}|)^6}$, then $|\lambda_{k+1}| \leq \|M_k\|_2 \leq \frac{|\lambda_{k+1}|}{1 - \delta_\times}$.*

(c) *Rayleigh quotient guarantee: if $\varepsilon_{\mathsf{pca}} \leq \frac{\delta_\times^5}{2^{28}k^3(|\lambda_1|/|\lambda_{k+1}|)^6}$, then $(1 - \delta_\times)|\lambda_k| \leq |v_k^\top Mv_k| \leq \frac{1}{1 - \delta_\times}|\lambda_k|$.*

(d) *Schatten-$q$ norm guarantee: for every $q \geq 1$, if $\varepsilon_{\mathsf{pca}} \leq \frac{\delta_\times^5}{2^{28}k^3d^{4/q}(|\lambda_1|/|\lambda_{k+1}|)^6}$, then*

$$\|M_k\|_{S_q} \leq \frac{(1 + \delta_\times)^2}{(1 - \delta_\times)^2}\Big(\sum_{i=k+1}^d \lambda_i^q\Big)^{1/q} = \frac{(1 + \delta_\times)^2}{(1 - \delta_\times)^2}\min_{V \in \mathbb{R}^{d \times k}, V^\top V = I}\big\{\|(I - VV^\top)M(I - VV^\top)\|_{S_q}\big\} .$$

---

*Proof of Theorem 4.1.* Let $V_s = (v_1, \ldots, v_s)$, so we can write

$$M_s = (I - V_sV_s^\top)M(I - V_sV_s^\top) = (I - v_sv_s^\top)M_{s-1}(I - v_sv_s^\top) .$$

We first claim that $\|M_{s-1}\|_2 \geq |\lambda_s|$ for every $s = 1, \ldots, k$. This can be proved by the Cauchy interlacing theorem. Indeed, $M_{s-1}^2 = (I - V_{s-1}V_{s-1}^\top)M^2(I - V_{s-1}V_{s-1}^\top)$ is a projection of $M^2$ into a $d - s + 1$ dimensional space, and therefore its largest eigenvalue $\|M_{s-1}^2\|_2$ should be at least as large as $|\lambda_s|^2$, the $s$-th largest eigenvalue of $M^2$. In other words, we have shown $\|M_{s-1}\|_2 \geq |\lambda_s|$.

(a) Define $\widehat{\lambda} = \|M_{k-1}\|_2 \geq |\lambda_k|$.

Note that all column vectors in $V_s$ are automatically eigenvectors of $M_s$ with eigenvalues zero. For analysis purpose only, let $W_s$ be the column matrix of eigenvectors in $V_s^\perp$ of $M_s$ that have *absolute* eigenvalues in the range $[0, (1 - \delta_\times + \tau_s)\widehat{\lambda}]$, where $\tau_s \stackrel{\text{def}}{=} \frac{s}{2k}\delta_\times$. We now show that for every $s = 0, 1, \ldots, k$, there exists a matrix $Q_s$ such that $\|U - W_s Q_s\|_2$ is small and $\|Q_s\|_2 \leq 1$. We will do this by induction.

In the base case: since $\tau_0 = 0$, we have $W_0 = U$ by the definition of $U$. We can therefore define $Q_0$ to be the identity matrix.

For every $s = 0, 1, \ldots, k-1$, suppose there exists a matrix $Q_s$ with $\|Q_s\|_2 \leq 1$ that satisfies $\|U - W_s Q_s\|_2 \leq \eta_s$ for some $\eta_s > 0$, we construct $Q_{s+1}$ as follows.

First we observe that `AppxPCA`$^\pm$ outputs a vector $v'_{s+1}$ satisfying $\|v'^\top_{s+1}W_s\|_2^2 \leq \varepsilon_{\mathsf{pca}}$ and $\|v'^\top_{s+1}V_s\|_2^2 \leq \varepsilon_{\mathsf{pca}}$ with probability at least $1 - p/k$. This follows from Theorem 3.1 (using $M = M_s$) because $[0, (1 - \delta_\times + \tau_s)\widehat{\lambda}] \subseteq [0, (1 - \delta_\times/2)\widehat{\lambda}]$, together with the fact that $\|M_s\|_2 \geq \|M_{k-1}\|_2 \geq \widehat{\lambda}$. Now, since $v_{s+1}$ is the projection of $v'_{s+1}$ into $V_s^\perp$, we have

$$\|v_{s+1}^\top W_s\|_2^2 \leq \frac{\|v'^\top_{s+1}W_s\|_2^2}{\|(I - V_s V_s^\top)v'_{s+1}\|_2^2} = \frac{\|v'^\top_{s+1}W_s\|_2^2}{1 - \|V_s^\top v'_{s+1}\|_2^2} \leq \frac{\varepsilon_{\mathsf{pca}}}{1 - \varepsilon_{\mathsf{pca}}} < 1.5\varepsilon_{\mathsf{pca}} \ . \tag{E.1}$$

Next we apply Lemma C.4 with $M = M_s$, $M' = M_{s+1}$, $U = W_s$, $V = W_{s+1}$, $v = v_{s+1}$, $\mu = (1 - \delta_\times + \tau_s)\widehat{\lambda}$, and $\tau = (\tau_{s+1} - \tau_s)\widehat{\lambda}$. We obtain a matrix $\widetilde{Q}_s$, $\|\widetilde{Q}_s\|_2 \leq 1$ such that[10]

$$\|W_s - W_{s+1}\widetilde{Q}_s\|_2 \leq \sqrt{\frac{169(\lambda_1/\widehat{\lambda})^2 \cdot 1.5\varepsilon_{\mathsf{pca}}}{(\tau_{s+1} - \tau_s)^2} + \varepsilon_{\mathsf{pca}}} < \frac{32\lambda_1 k\sqrt{\varepsilon_{\mathsf{pca}}}}{\lambda_k \delta_\times} \ ,$$

and this implies that

$$\|W_{s+1}\widetilde{Q}_s Q_s - U\|_2 \leq \|W_{s+1}\widetilde{Q}_s Q_s - W_s Q_s\|_2 + \|W_s Q_s - U\|_2 \leq \eta_s + \frac{32\lambda_1 k\sqrt{\varepsilon_{\mathsf{pca}}}}{\lambda_k \delta_\times} \ .$$

Let $Q_{s+1} = \widetilde{Q}_s Q_s$ we know that $\|Q_{s+1}\|_2 \leq 1$ and

$$\|W_{s+1}Q_{s+1} - U\|_2 \leq \eta_{s+1} \stackrel{\text{def}}{=} \eta_s + \frac{32\lambda_1 k\sqrt{\varepsilon_{\mathsf{pca}}}}{\lambda_k \delta_\times} \ .$$

Therefore, after $k$-iterations of `LazyEV`, we obtain:

$$\|W_k Q_k - U\|_2 \leq \eta_k = \frac{32\lambda_1 k^2 \sqrt{\varepsilon_{\mathsf{pca}}}}{\lambda_k \delta_\times}$$

Multiply $U^\top$ from the left, we obtain $\|U^\top W_k Q_k - I\|_2 \leq \eta_k$. Since $\|Q_k\|_2 \leq 1$, we must have $\sigma_{\min}(U^\top W_k) \geq 1 - \eta_k$ (here $\sigma_{\min}$ denotes the smallest singular value). Therefore,

$$U^\top W_k W_k^\top U \succeq (1 - \eta_k)^2 I \ .$$

Since $V_k$ and $W_k$ are orthogonal of each other, we have

$$U^\top V_k V_k^\top U \preceq U^\top (I - W_k W_k^\top)U \preceq I - (1 - \eta_k)^2 I \preceq 2\eta_k I$$

---

[10] Technically speaking, to apply Lemma C.4 we need $U = W_s$ to consist of all eigenvectors of $M_s$ with absolute eigenvalues $\leq \mu$. However, we only defined $W_s$ to be such eigenvectors that are *also* orthogonal to $V_s$. It is straightforward to verify that the same result of Lemma C.4 remains true because $v_{s+1}$ is orthogonal to $V_s$.

Therefore,

$$\|V_k^\top U\|_2 \le 8 \frac{(|\lambda_1|/|\lambda_k|)^{1/2} k \varepsilon_{\mathsf{pca}}^{1/4}}{\delta_\times^{1/2}} \le \varepsilon \ .$$

(b) The statement is obvious when $k = 0$. For every $k \ge 1$, the lower bound is obvious. We prove the upper bound by contradiction. Suppose that $\|M_k\|_2 > \frac{|\lambda_{k+1}|}{1-\delta_\times}$. Then, since $\|M_{k-1}\|_2 \ge \|M_k\|_2$ and therefore $|\lambda_{k+1}|, \ldots, |\lambda_d| < (1-\delta_\times)\|M_{k-1}\|_2$, we can apply Theorem 4.1.a of the current $k$ to deduce that $\|V_k^\top U_{>k}\|_2 \le \varepsilon$ where $U_{>k} \stackrel{\text{def}}{=} (u_{k+1}, \ldots, u_d)$. We now apply Lemma C.2 with $V_s = V_k$ and $U = U_{>k}$, we obtain a matrix $Q_k \in \mathbb{R}^{d \times k}$ whose columns are spanned by $u_1, \ldots, u_k$ and satisfy

$$\left\| \left(I - Q_k Q_k^\top\right) M \left(I - Q_k Q_k^\top\right) - \left(I - V_k V_k^\top\right) M \left(I - V_k V_k^\top\right) \right\|_2 < 16|\lambda_1|\varepsilon \ .$$

However, our assumption says that the second matrix $\left(I - V_k V_k^\top\right) M \left(I - V_k V_k^\top\right)$ has spectral norm at least $|\lambda_{k+1}|/(1-\delta_\times)$, but we know that $\left(I - Q_k Q_k^\top\right) M \left(I - Q_k Q_k^\top\right)$ has spectral norm exactly $|\lambda_{k+1}|$ due to the definition of $Q_k$. Therefore, we must have $\frac{|\lambda_{k+1}|}{1-\delta_\times} - |\lambda_{k+1}| \le 16|\lambda_1|\varepsilon$ due to triangle inequality.

In other words, by selecting $\varepsilon$ in Theorem 4.1.a to satisfy $\varepsilon \le \frac{\delta_\times}{16|\lambda_1|/|\lambda_{k+1}|}$ (which is satisfied by our assumption on $\varepsilon_{\mathsf{pca}}$), we get a contradiction so can conclude that $\|M_k\|_2 \le \frac{|\lambda_{k+1}|}{1-\delta_\times}$.

(c) We compute that

$$
\begin{aligned}
|v_k^\top M v_k| = |v_k^\top M_{k-1} v_k| &\overset{\text{①}}{\ge} \frac{|v_k'^\top M_{k-1} v_k'|}{\|(I - V_{k-1}V_{k-1}^\top)v_k'\|_2^2} \overset{\text{②}}{\ge} \frac{|v_k'^\top M_{k-1} v_k'|}{(1 - \sqrt{\varepsilon_{\mathsf{pca}}})^2} \\
&\overset{\text{③}}{\ge} \frac{1 - \varepsilon_{\mathsf{pca}}}{(1 - \sqrt{\varepsilon_{\mathsf{pca}}})^2}(1 - \delta_\times/2)\|M_{k-1}\|_2 \ge (1 - \delta_\times)\|M_{k-1}\|_2 \ . \quad (\text{E.2})
\end{aligned}
$$

Above, ① is because $v_k$ is the projection of $v_k'$ into $V_{k-1}^\perp$, ② is because $\|V_{k-1}^\top v_k'\|_2^2 \le \varepsilon_{\mathsf{pca}}$ following the same reason as (E.1), and ③ is owing to Theorem 3.1. Next, since $\|M_{k-1}\|_2 \ge |\lambda_k|$, we automatically have $|v_k^\top M v_k| \ge (1 - \delta_\times)|\lambda_k|$. On the other hand, $|v_k^\top M v_k| = |v_k^\top M_{k-1} v_k| \le \|M_{k-1}\|_2 \le \frac{|\lambda_k|}{1-\delta_\times}$ where the last inequality is owing to Theorem 4.1.b.

(d) Since $\|V_k^\top U\|_2 \le \varepsilon_c \stackrel{\text{def}}{=} 8\frac{(|\lambda_1|/|\lambda_k|)^{1/2} k \varepsilon_{\mathsf{pca}}^{1/4}}{\delta_\times^{1/2}}$ from the analysis of Theorem 4.1.a, we can apply Lemma C.2 to obtain a (column) orthogonal matrix $Q_k \in \mathbb{R}^{d \times k}$ such that

$$\|M_k' - M_k\|_2 \le 16|\lambda_1|\varepsilon_c, \qquad \text{where } M_k' \stackrel{\text{def}}{=} (I - Q_k Q_k^\top) M (I - Q_k Q_k^\top) \qquad (\text{E.3})$$

Suppose $U = (u_{d-p+1}, \ldots, u_d)$ is of dimension $d \times p$, that is, there are exactly $p$ eigenvalues of $M$ whose absolute value is $\le (1 - \delta_\times)\|M_{k-1}\|_2$. Then, the definition of $Q_k$ in Lemma C.2 tells us $U^\top Q_k = 0$ so $M_k'$ agrees with $M$ on all the eigenvalues and eigenvectors $\{(\lambda_j, u_j)\}_{j=d-p+1}^d$ because an index $j$ satisfies $|\lambda_j| \le (1 - \delta_\times)\|M_{k-1}\|_2$ if and only if $j \in \{d-p+1, d-p+2, \ldots, d\}$.

Denote by $\mu_1, \ldots, \mu_{d-k}$ the eigenvalues of $M_k'$ excluding the $k$ zero eigenvalues in subspace $Q_k$, and assume without loss of generality that $\{\mu_1, \ldots, \mu_p\} = \{\lambda_{d-p+1}, \ldots, \lambda_d\}$. Then,

$$
\begin{aligned}
\|M_k'\|_{S_q}^q &= \sum_{i=1}^{d-k} |\mu_i|^q = \sum_{i=1}^{p} |\mu_i|^q + \sum_{i=p+1}^{d-k} |\mu_i|^q = \sum_{i=d-p+1}^{d} |\lambda_i|^q + \sum_{i=p+1}^{d-k} |\mu_i|^q \\
&\overset{\text{①}}{\le} \sum_{i=d-p+1}^{d} |\lambda_i|^q + (d-k-p)\|M_k'\|_2^q \overset{\text{②}}{\le} \sum_{i=d-p+1}^{d} |\lambda_i|^q + (d-k-p)(\|M_k\|_2 + 16|\lambda_1|\varepsilon_c)^q
\end{aligned}
$$

$$\overset{\text{③}}{\leq} \sum_{i=d-p+1}^{d} |\lambda_i|^q + (d-k-p)\Big(\frac{|\lambda_{k+1}|}{(1-\delta_\times)} + 16|\lambda_1|\varepsilon_c\Big)^q$$

Above, ① is because each $|\mu_i|$ is no greater than $\|M'_k\|_2$, and ② is owing to (E.3), and ③ is because of Theorem 4.1.b. Suppose we choose $\varepsilon_c$ so that $\varepsilon_c \leq \frac{|\lambda_{k+1}|\delta_\times}{16\lambda_1}$ (and this is indeed satisfied by our assumption on $\varepsilon_{\texttt{pca}}$), then we can continue and write

$$\|M'_k\|_{S_q}^q \leq \sum_{i=d-p+1}^{d} |\lambda_i|^q + (d-k-p)\frac{(1+\delta_\times)^q}{(1-\delta_\times)^q}|\lambda_{k+1}|^q$$

$$\overset{\text{④}}{\leq} \sum_{i=d-p+1}^{d} |\lambda_i|^q + \frac{(1+\delta_\times)^q}{(1-\delta_\times)^{2q}}\sum_{i=k+1}^{d-p} |\lambda_i|^q \leq \frac{(1+\delta_\times)^q}{(1-\delta_\times)^{2q}}\sum_{i=k+1}^{d} |\lambda_i|^q \ .$$

Above, ④ is because for each eigenvalue $\lambda_i$ where $i \in \{k+1, k+2, \ldots, d-p\}$, we have $|\lambda_i| > (1-\delta_\times)\|M_{k-1}\|_2 \geq (1-\delta_\times)|\lambda_k| \geq (1-\delta_\times)|\lambda_{k+1}|$. Finally, using (E.3) again we have

$$\|M_k\|_{S_q} \leq \|M'_k\|_{S_q} + \|M_k - M'_k\|_{S_q} \leq \|M'_k\|_{S_q} + d^{1/p}\|M_k - M'_k\|_2$$

$$\leq \frac{1+\delta_\times}{(1-\delta_\times)^2}\Big(\sum_{i=k+1}^{d} |\lambda_i|^q\Big)^{1/q} + 16d^{1/p}|\lambda_1|\varepsilon_c$$

As long as $\varepsilon_c \leq \frac{\delta_\times|\lambda_{k+1}|}{16d^{1/p}\lambda_1}$, we have

$$\|M_k\|_{S_q} \leq \frac{(1+\delta_\times)^2}{(1-\delta_\times)^2}\Big(\sum_{i=k+1}^{d} |\lambda_i|^q\Big)^{1/q}$$

as desired. Finally, we note that $\varepsilon_c \leq \frac{\delta_\times\lambda_{k+1}}{16d^{1/p}\lambda_1}$ is satisfied with our assumption on $\varepsilon_{\texttt{pca}}$, and note that $\min_{V\in\mathbb{R}^{d\times k}, V^\top V=I}\big\{\|(I-VV^\top)M(I-VV^\top)\|_{S_q}\big\} = \big(\sum_{i=k+1}^{d} |\lambda_i|^q\big)^{1/q}$ which follows easily from Cauchy interlacing theorem. $\qquad\square$

## E.2 Proofs of Theorems 4.3 and 4.4

*Proof of Theorem 4.3.* Define $V_k = B^{1/2}\mathbb{V}_k = \texttt{LazyEV}(\cdots)$ to be the direct output of $\texttt{LazyEV}$. The column orthogonality of $V_k$ implies $\mathbb{V}_k^\top B\mathbb{V}_k = I$.

It is clear from the definition of generalized eigenvectors that $B^{1/2}u_1, \ldots, B^{1/2}u_d$ are eigenvectors of $M \overset{\text{def}}{=} B^{-1/2}AB^{-1/2}$ with eigenvalues $\lambda_1, \ldots, \lambda_d$. Applying Theorem 4.1.a, we have: $\|V_k^\top U\|_2 \leq \varepsilon$ where $U = (B^{1/2}u_j, \ldots, B^{1/2}u_d)$ is a (column) orthonormal matrix and $j$ is the smallest index satisfying $|\lambda_j| \leq (1-\delta_\times)\|M_{k-1}\|_2$. Since it satisfies $\|M_{k-1}\|_2 \geq |\lambda_k|$, we have

$$|\lambda_{k+1}| = |\lambda_k|(1-\texttt{gap}) = |\lambda_k|(1-\delta_\times) \leq (1-\delta_\times)\|M_{k-1}\|_2 \ .$$

Therefore, $j$ must be equal to $k+1$ according to its definition, so we have $U = B^{1/2}\mathbb{W}$. This implies $\|\mathbb{V}_k^\top B\mathbb{W}\|_2 = \|V_k^\top U\|_2 \leq \varepsilon$.

The running time statement comes directly from Theorem 4.2 by putting in the parameters. $\qquad\square$

*Proof of Theorem 4.4.* Define $V_k = B^{1/2}\mathbb{V}_k = \texttt{LazyEV}(\cdots)$ to be the direct output of $\texttt{LazyEV}$. The column orthogonality of $V_k$ implies $\mathbb{V}_k^\top B\mathbb{V}_k = I$.

It is clear from the definition of generalized eigenvectors that $B^{1/2}u_1, \ldots, B^{1/2}u_d$ are eigenvectors of $M \overset{\text{def}}{=} B^{-1/2}AB^{-1/2}$ with eigenvalues $\lambda_1, \ldots, \lambda_d$. Applying Theorem 4.1.b, we have: $\big\|(I - V_kV_k^\top)B^{-1/2}AB^{-1/2}(I-V_kV_k^\top)\big\|_2 \leq \frac{|\lambda_{k+1}|}{1-\varepsilon}$. Next, for every vector $\mathbb{w} \in \mathbb{R}^d$ that is $B$-orthogonal to

$\mathbb{V}_k$, that is, $\mathtt{w}^\top B \mathbb{V}_k = 0$, we can define $w \overset{\text{def}}{=} B^{1/2}\mathtt{w}$ and we know $w$ is orthogonal to $V_k$. We can apply the above spectral upper bound and get

$$\mathtt{w}^\top A\mathtt{w} = w^\top B^{-1/2}AB^{-1/2}w = w^\top(I - V_k V_k^\top)B^{-1/2}AB^{-1/2}(I - V_k V_k^\top)w$$
$$\leq \|w\|_2^2 \cdot \frac{|\lambda_{k+1}|}{1-\varepsilon} = \mathtt{w}^\top B\mathtt{w} \cdot \frac{|\lambda_{k+1}|}{1-\varepsilon}$$

as desired. At the same time, denoting by $v_s = B^{1/2}\mathtt{v}_s$, Theorem 4.1.c implies that

$$\forall s \in [k]: \quad \left|\mathtt{v}_s^\top A\mathtt{v}_s\right| = \left|v_s^\top M v_s\right| \in \left[(1-\varepsilon)|\lambda_s|, \frac{|\lambda_s|}{1-\varepsilon}\right] \ .$$

The running time statement comes directly from Theorem 4.2 by putting in the parameters. $\qquad\square$

# F  Proof for Section 6: CCA Theorems

## F.1  The Main Convergence Theorem

Since $\mathtt{LazyCCA}$ only admits minor changes on top of $\mathtt{LazyEV}$, the next theorem is an almost identical copy of Theorem 4.1. To make this paper concise, instead of reproving Theorem F.1 line by line, we here only sketch the main changes needed in the new proof.

---

**Theorem F.1.** *Let $M = B^{-1/2}AB^{-1/2} \in \mathbb{R}^{d \times d}$ be a symmetric matrix where $A$ and $B$ are matrices coming from a CCA instance using Lemma 2.3. Suppose $M$ has eigenvalues $\lambda_1, \dots, \lambda_d \in [-1, 1]$ and corresponding eigenvectors $u_1, \dots, u_d$. Suppose without loss of generality that $|\lambda_1| \geq \cdots \geq |\lambda_d|$.*

*For every $k \in [d]$, $\delta_\times, p \in (0, 1)$, there exists some $\varepsilon_{\mathtt{pca}} \leq O\big(\mathtt{poly}(\delta_\times, \frac{|\lambda_1|}{|\lambda_{k+1}|}, \frac{1}{d})\big)$ such that $\mathtt{LazyCCA}$ outputs a (column) orthonormal matrix $V_k = (v_1, \dots, v_{2k}) \in \mathbb{R}^{d \times 2k}$ which, with probability at least $1-p$, satisfies all of the following properties. (Denote by $M_s = (I - V_s V_s^\top)M(I - V_s V_s^\top)$.)*

*(a) Correlation guarantee: $\|V_k^\top U\|_2 \leq \varepsilon$,*
    *where $U = (u_j, \dots, u_d)$ and $j$ is the smallest index satisfying $|\lambda_j| \leq (1-\delta_\times)\|M_{k-1}\|_2$.*

*(b) Spectral norm guarantee: $|\lambda_{2k+1}| \leq \|M_k\|_2 \leq \frac{|\lambda_{2k+1}|}{1-\delta_\times}$.*

*(c) Rayleigh quotient guarantee: $(1-\delta_\times)|\lambda_{2k}| \leq |v_{2k-1}^\top M v_{2k-1}| = |v_{2k}^\top M v_{2k}| \leq \frac{1}{1-\delta_\times}|\lambda_{2k}|$.*

*(d) Schatten-$q$ norm guarantee: for every $q \geq 1$,*

$$\|M_k\|_{S_q} \leq \frac{(1+\delta_\times)^2}{(1-\delta_\times)^2}\Big(\sum_{i=2k+1}^d \lambda_i^q\Big)^{1/q} = \frac{(1+\delta_\times)^2}{(1-\delta_\times)^2} \min_{V \in \mathbb{R}^{d \times 2k}, V^\top V = I}\big\{\|(I-VV^\top)M(I-VV^\top)\|_{S_q}\big\}.$$

---

*Proof sketch of Theorem F.1.* Recall that when a vector $v_s \in \mathbb{R}^d$ is obtained in iteration $s$ of $\mathtt{LazyEV}$, the proof of Theorem 4.1 suggest that the following two properties hold

$$\left|v_s^\top M_{s-1} v_s\right| \geq (1-\delta_\times)\|M_{s-1}\|_2 \quad \text{and} \quad \left\|v_s^\top W_{s-1}\right\|_2^2 \leq 1.5\varepsilon_{\mathtt{pca}} \ . \tag{F.1}$$

(The first property is shown in (E.2), and the second property is shown in (E.1). Recall that $W_{s-1}$ is the column orthonormal matrix containing all eigenvectors of $M_{s-1}$ whose absolute eigenvalues are below some threshold.) Then, the proof of Theorem 4.1 proceeds by heavily relying on (F.1).

In our $\mathtt{LazyCCA}$, after obtaining this same vector $v_s$, we write it as $v_s = (\xi_s', \zeta_v')$ and perform block-scaling $\xi_s = \xi_s'/(\sqrt{2}\|\xi_s'\|)$ and $\zeta_s = \zeta_s'/(\sqrt{2}\|\zeta_s'\|)$, see Line 7 of $\mathtt{LazyCCA}$. Therefore, in order for the same proof of Theorem 4.1 to hold, we need to show that this new vector $(\xi_s, \zeta_s)$ satisfies the same properties up to constants:

$$\left| \left( \begin{array}{c} \xi_s \\ \zeta_s \end{array} \right)^\top M_{s-1} \left( \begin{array}{c} \xi_s \\ \zeta_s \end{array} \right) \right| \geq (1-\delta_\times)\|M_{s-1}\|_2 \quad \text{and} \quad \left\| \left( \begin{array}{c} \xi_s \\ \zeta_s \end{array} \right)^\top W_{s-1} \right\|_2^2 \leq 12.5\varepsilon_{\mathsf{pca}} \ . \tag{F.2}$$

Suppose $V_{s-1} = \left( \begin{array}{ccc} \pm\xi_1 & \cdots & \pm\xi_{s-1} \\ \zeta_1 & \cdots & \zeta_{s-1} \end{array} \right)$. Since $v_s$ is orthogonal to all vectors in the column span of $V_{s-1}$ according to Line 5 of `LazyCCA`, we automatically have $\xi_s^\top \xi_i = 0$ and $\zeta_s^\top \zeta_i = 0$ for all $i \in [s-1]$. We also have $\|\xi_s\|^2 + \|\zeta_s\|^2 = 1/2 + 1/2 = 1$ so the new vector $(\xi_s, \zeta_s)$ has Euclidean norm 1.

As for the first property in (F.2), we observe that the new vector $(\xi_s, \zeta_s)$ enjoys an (absolute) Rayleigh quotient value that is no worse than the original $v_s = (\xi_s', \zeta_s')$. This is so because (without loss of generality we consider $v_s^\top M_{s-1} v_s > 0$):

$$\left( \begin{array}{c} \xi_s \\ \zeta_s \end{array} \right)^\top M_{s-1} \left( \begin{array}{c} \xi_s \\ \zeta_s \end{array} \right) \overset{①}{=} \left( \begin{array}{c} \xi_s \\ \zeta_s \end{array} \right)^\top M \left( \begin{array}{c} \xi_s \\ \zeta_s \end{array} \right) \overset{②}{=} 2\xi_s^\top \big(S_{xx}^{-1/2} S_{xy} S_{yy}^{-1/2}\big) \zeta_s$$

$$\overset{③}{=} \frac{1}{\|\xi_s'\|\|\zeta_s'\|} \cdot \xi_s'^\top \big(S_{xx}^{-1/2} S_{xy} S_{yy}^{-1/2}\big) \zeta_s' \overset{④}{=} \frac{1}{2\|\xi_s'\|\|\zeta_s'\|} \cdot v_s^\top M v_s \overset{⑤}{=} \frac{1}{2\|\xi_s'\|\|\zeta_s'\|} \cdot v_s^\top M_{s-1} v_s \overset{⑥}{\geq} v_s^\top M_{s-1} v_s \ . \tag{F.3}$$

Above, ① is because $(\xi_s, \zeta_s)$ is orthogonal to $V_{s-1}$; ② is by the definition of $M = B^{-1/2}AB^{-1/2}$ as well as the definition of $A$ and $B$; ③ is by the definitions of $\xi_s$ and $\zeta_s$; ④ is by $v_s = (\xi_s', \zeta_s')$ and again by the definition of $M$; ⑤ follows from the fact that $v_s$ is orthogonal to $V_{s-1}$; and ⑥ follows from AM-GM together with the fact that $\|\xi_s'\|^2 + \|\zeta_s'\|^2 = \|v_s\|^2 = 1$. This finishes proving the first property in (F.2) because the original vector $v_s$ satisfies $|v_s^\top M_{s-1} v_s| \geq (1-\delta_\times)\|M_{s-1}\|_2$ according to (F.1).

We make an additional observation here: $\|\xi_s'\|^2$ and $\|\zeta_s'\|^2$ must be in the range $[0.06, 0.94]$ before scaling. Indeed, suppose for instance $\|\xi_s'\|^2 = c$ for some $c \in [0,1]$. Then, it satisfies $2\|\xi_s'\|\|\zeta_s'\| = 2\sqrt{c(1-c)}$ and therefore (F.3) becomes $\left( \begin{array}{c} \xi_s \\ \zeta_s \end{array} \right)^\top M_{s-1} \left( \begin{array}{c} \xi_s \\ \zeta_s \end{array} \right) \geq \frac{1-\delta_\times}{2\sqrt{c(1-c)}}\|M_{s-1}\|_2$, meaning that $\frac{1-\delta_\times}{2\sqrt{c(1-c)}} \leq 1$. If $\delta_\times \leq 1/2$, this implies $c - 1/2 \in [-\sqrt{3}/4, \sqrt{3}/4]$ and thus $c \in [0.06, 0.94]$.

As for the second property in (F.2), for every matrix $W_{s-1}$ that is in the proof of Theorem 4.1.a, its columns are all eigenvectors of $M_{s-1}$ whose absolute eigenvalues are below some threshold, so must consist of only symmetric vectors in this CCA setting: that is, $W_{s-1} = \left( \begin{array}{ccc} \pm a_1 & \cdots & \pm a_t \\ b_1 & \cdots & b_t \end{array} \right).$[11] According to (F.1) we already know $\|v_s^\top W_{s-1}\|^2 \leq 1.5\varepsilon_{\mathsf{pca}}$, which implies

$$\|v_s^\top W_{s-1}\|^2 = \sum_{i=1}^t (\xi_s'^\top a_i + \zeta_s'^\top b_i)^2 + (\xi_s'^\top a_i - \zeta_s'^\top b_i)^2 = 2\|\xi_s'^\top (a_1, \ldots, a_t)\|^2 + 2\|\zeta_s'^\top (b_1, \ldots, b_t)\|^2 \leq 1.5\varepsilon_{\mathsf{pca}} \ .$$

Now we can compute

$$\left\| \left( \begin{array}{c} \xi_s \\ \zeta_s \end{array} \right)^\top W_{s-1} \right\|_2^2 = 2\|\xi_s^\top (a_1, \ldots, a_t)\|^2 + 2\|\zeta_s^\top (b_1, \ldots, b_t)\|^2$$

$$\leq \frac{0.5}{0.06} \cdot \big( 2\|\xi_s'^\top (a_1, \ldots, a_t)\|^2 + 2\|\zeta_s'^\top (b_1, \ldots, b_t)\|^2 \big) \leq 12.5\varepsilon_{\mathsf{pca}} \ ,$$

---

[11] This is because matrix $M_s$ is always of the form $D^{-1/2}CD^{1/2}$ where $D = \mathrm{diag}\{D_1, D_2\}$ is block diagonal and positive definite, while $C = [[0, C_1]; [C_1^\top, 0]]$ has only zero on its two block diagonal locations. The same proof of Lemma 2.3 shows that the eigenvectors of $D^{-1/2}CD^{1/2}$ must be symmetric. In fact, to be precise, $W_s$ may also contain some eigenvectors corresponding to zero eigenvalues. However, adding them will make our notations heavier, so we refrain from doing that in this sketched proof.

where the first inequality is because $\|\xi'_s\|^2, \|\zeta'_s\|^2 \geq 0.06$.

This finishes proving the two properties in (F.2), so Theorem F.1 holds after plugging the rest of the proof of Theorem 4.1 in but changing constants slightly. $\square$

## F.2  Fast Implementation of LazyCCA: Stochastic

> **Theorem F.2** (stochastic running time of `LazyCCA`). *Let $X \in \mathbb{R}^{n \times d_x}, Y \in \mathbb{R}^{n \times d_y}$ be two matrices, and define $A$ and $B$ according to Lemma 2.3. Suppose $M = B^{-1/2}AB^{-1/2}$, and `RanInit`$(d)$ is the random vector generator defined in Proposition 3.3, and we want to compute matrix $\mathbb{V} \leftarrow B^{-1/2}$`LazyCCA`$(\mathcal{A}, M, k, \delta_\times, \varepsilon_{\mathsf{pca}}, p)$. Then, this procedure can be implemented to run in time*
>
> - $\widetilde{O}\left(\frac{k \operatorname{nnz}(B) + k^2 d + k \Upsilon}{\sqrt{\delta_\times}}\right)$ *where $\Upsilon$ is the time needed to compute $B^{-1}Aw$ for a vector $w$ to an accuracy $\varepsilon$ where $\log(1/\varepsilon) = \widetilde{O}(1)$, or*
>
> - $\widetilde{O}\left(\frac{k \operatorname{nnz}(X,Y)\left(1 + \sqrt{\kappa'/n}\right) + k^2 d}{\sqrt{\delta_\times}}\right)$ *if we simply use `Katyusha` to compute $B^{-1}Aw$.*
>
> *Above, $\kappa = \frac{\lambda_{\max}(B)}{\lambda_{\min}(B)} = \frac{\max\{\lambda_{\max}(S_{xx}), \lambda_{\max}(S_{yy})\}}{\min\{\lambda_{\min}(S_{xx}), \lambda_{\min}(S_{yy})\}}$, and $\kappa' = \frac{2 \max_{i \in [n]} \{\|X_i\|^2, \|Y_i\|^2\}}{\lambda_{\min}(B)} \in [\kappa, 2n\kappa]$.*

*Proof.* The proof of the first item is almost identical to the proof of Theorem 4.2 so ignored here. As for the second item, it follows from the first item together with the running time of `Katyusha` in Lemma 2.6. $\square$

## F.3  Fast Implementation of LazyCCA: Doubly Stochastic

> **Theorem F.3** (doubly-stochastic running time of `LazyCCA`). *Let $X \in \mathbb{R}^{n \times d_x}, Y \in \mathbb{R}^{n \times d_y}$ be two matrices, and define $A$ and $B$ according to Lemma 2.3. Suppose $M = B^{-1/2}AB^{-1/2}$, and `RanInit`$(d)$ is the random vector generator defined in Proposition 3.3, and we want to compute $\mathbb{V} \leftarrow B^{-1/2}$`LazyCCA`$(\mathcal{A}, M, k, \delta_\times, \varepsilon_{\mathsf{pca}}, p)$. Then, this procedure can be implemented to run in time*
>
> - $\widetilde{O}\left(\operatorname{nnz}(X,Y) \cdot \left(1 + \frac{\sqrt{\kappa'/n^{1/4}}}{\sqrt{\delta_\times \sigma_1}}\right)\right)$ *if $k = 1$ and we use accelerated SVRG as the method $\mathcal{A}$;*
>
> - $\widetilde{O}\left(k \operatorname{nnz}(X,Y) \cdot \left(1 + \frac{\sqrt{\kappa'}}{\sqrt{\delta_\times \sigma_k(\operatorname{nnz}(X,Y)/kd)^{1/4}}}\right)\right)$ *if $\operatorname{nnz}(X,Y) \geq kd$ and we use accelerated SVRG as the method $\mathcal{A}$.*
>
> *Above, $\kappa = \frac{\lambda_{\max}(B)}{\lambda_{\min}(B)} = \frac{\max\{\lambda_{\max}(S_{xx}), \lambda_{\max}(S_{yy})\}}{\min\{\lambda_{\min}(S_{xx}), \lambda_{\min}(S_{yy})\}}$, and $\kappa' = \frac{2 \max_{i \in [n]} \{\|X_i\|^2, \|Y_i\|^2\}}{\lambda_{\min}(B)} \in [\kappa, 2n\kappa]$.*

To prove Theorem F.3, let us recall from the proof of Theorem 4.2 (see Section 5.2) that it suffices to bound the time needed to compute $\mathbb{N}^{-1}w$ where

$$\mathbb{N} \overset{\text{def}}{=} B^{-1/2}(\lambda I - M_s)B^{1/2} = \lambda I - (I - \mathbb{V}_s \mathbb{V}_s^\top B)B^{-1}A(I - \mathbb{V}_s \mathbb{V}_s^\top B) ~,$$

and we only need to compute it poly-logarithmic number of times for each $s = 1, \ldots, k$.

Observe now that $\mathbb{N}^{-1}w = (B\mathbb{N})^{-1}Bw$, so it suffices to bound the time needed to compute $(B\mathbb{N})^{-1}$ applied to a vector, and we have (dropping the subscript of $\mathbb{V}$ for cleanness)

$$B\mathbb{N} = \lambda B - (I - B\mathbb{V}\mathbb{V}^\top)A(I - \mathbb{V}\mathbb{V}^\top B) ~.$$

### F.3.1  Special Case of 1-CCA

When $k = 1$, the matrix $\mathbb{V}$ is empty and we want to compute $(\lambda B - A)^{-1}w$ for some vector $w$. This is equivalent to minimizing $f(z) \overset{\text{def}}{=} \frac{1}{2}z^\top(\lambda B - A)z - w^\top z$. We analyze its running time as follows, and a similar analysis has also appeared in [29, Section 3.2.3].

Using the definition of $A$ and $B$, we can write $f(z) = \frac{1}{n}\sum_{i=1}^{n} f_i(z)$ where

$$f_i(z) \stackrel{\text{def}}{=} \frac{1}{2}\lambda(\langle X_i, z_1\rangle^2 + \langle Y_i, z_2\rangle^2) - 2\langle X_i, z_1\rangle\langle Y_i, z_2\rangle + w^\top z\ ,$$

where we have denoted by $z = (z_1, z_2)$ for $z_1 \in \mathbb{R}^{d_x}$ and $z_2 \in \mathbb{R}^{d_y}$, and by $X_i, Y_i$ the $i$-th row vector of matrix $X$ and $Y$ respectively. The smoothness of $f_i(z)$ is given by (using the fact that $\lambda \leq 2$ which comes from the definition of $\texttt{AppxPCA}^{\pm}$)

$$\|\nabla^2 f_i(z)\|_2 \leq O(1) \cdot \max\{\|X_i\|^2, \|Y_i\|^2\}\ .$$

The strong convexity of $f(z)$ is given by

$$\nabla^2 f(z) = B^{1/2}(\lambda I - B^{-1/2}AB^{-1/2})B^{1/2} \succeq B^{1/2}(\frac{\delta_\times\sigma_1}{48}I)B^{1/2} \succeq \frac{\sigma_1\delta_\times\lambda_{\min}(B)}{48}\cdot I\ ,$$

where the first inequality is because $\lambda_{\min}(\lambda I - B^{-1/2}AB^{-1/2}) \geq \frac{\delta_\times\lambda^*}{48}$ from Theorem 3.1, as well as $\lambda^* \geq \sigma_1 = \lambda_{\max}(B^{-1/2}AB^{-1/2})$ from the definition of the algorithm.

Therefore, by applying the accelerated SVRG method on minimizing the sum-of-non-convex objective $f(z)$ [7, 24], we can find an $\widetilde{\varepsilon}$ minimizer of $f(z)$ in time

$$\widetilde{O}\Big(\texttt{nnz}(X,Y)\cdot\Big(1 + \frac{\sqrt{\max_i\{\|X_i\|^2, \|Y_i\|^2\}}}{\sqrt{\sigma_1\delta_\times\lambda_{\min}(B)}n^{1/4}}\Big)\cdot\log^2\frac{\|w\|}{\widetilde{\varepsilon}}\Big) = \widetilde{O}\Big(\texttt{nnz}(X,Y)\cdot\big(1 + \frac{\sqrt{\kappa'}}{\sqrt{\sigma_1\delta_\times}n^{1/4}}\big)\cdot\log^2\frac{\|w\|}{\widetilde{\varepsilon}}\Big)\ .$$

Finally, similar to the proof of Lemma 2.6, an $\widetilde{\varepsilon}$ minimizer of $f(z)$ implies $\|(\lambda B - A)^{-1}w - z\| \leq \varepsilon$ when $\varepsilon^2 = 96\widetilde{\varepsilon}/(\sigma_1\delta_\times\lambda_{\min}(B))$. Putting this back we obtain a final running time of $\widetilde{O}\Big(\texttt{nnz}(X,Y)\cdot\big(1 + \frac{\sqrt{\kappa'}}{\sqrt{\sigma_1\delta_\times}n^{1/4}}\big)\Big)$, ignoring log factors.

### F.3.2 Genearl Case of $k$-CCA

For the general case when $\mathbb{V}$ is not empty, one can carefully check that to compute

$$(\lambda B - (I - B\mathbb{V}\mathbb{V}^\top)A(I - \mathbb{V}\mathbb{V}^\top B))^{-1}w = (\lambda B - A + (B\mathbb{V}\mathbb{V}^\top A + A\mathbb{V}\mathbb{V}^\top B - B\mathbb{V}\mathbb{V}^\top A\mathbb{V}\mathbb{V}^\top B)^{-1}w$$

it suffices to minimize $f'(z) = \frac{1}{n}\sum_{i=1}^{n} f_i'(z)$ for

$$f_i'(z) \stackrel{\text{def}}{=} f_i(z) + \frac{1}{2}z^\top Qz \quad \text{where} \quad Q \stackrel{\text{def}}{=} B\mathbb{V}\mathbb{V}^\top A + A\mathbb{V}\mathbb{V}^\top B - B\mathbb{V}\mathbb{V}^\top A\mathbb{V}\mathbb{V}^\top B\ ,$$

and $f_i(z)$ is the same as defined in the $k = 1$ case of the previous subsection.

**Claim F.4.** $\|\nabla^2 f_i'(z)\|_2 \leq O(1)\cdot\max\{\|X_i\|^2, \|Y_i\|^2\}$ and $\nabla^2 f'(z) \succeq \frac{\sigma_k\delta_\times\lambda_{\min}(B)}{48}\cdot I$ .

*Proof.* Since

$$\|B\mathbb{V}\|_2 \leq \|B^{1/2}\|_2\|B^{-1/2}V\|_2 \leq \|B\|_2^{1/2}$$

$$\|A\mathbb{V}\|_2 \leq \|B^{1/2}\|_2\|B^{-1/2}AB^{-1/2}\|_2\|B^{1/2}\mathbb{V}\|_2 \leq \|B\|_2^{1/2}\sigma_1$$

$$\|\mathbb{V}^\top A\mathbb{V}^\top\|_2 \leq \|\mathbb{V}^\top B^{1/2}\|_2\|B^{-1/2}AB^{-1/2}\|_2\|B^{1/2}\mathbb{V}\|_2 \leq \sigma_1$$

we can easily compute that $\|Q\|_2 \leq O(1)\cdot\|B\|_2 \leq O(1)\cdot\text{Tr}(B) \leq O(1)\cdot\max_i\{\|X_i\|^2, \|Y_i\|^2\}$. Therefore, we have the same smoothness property as in the previous subsection:

$$\|\nabla^2 f_i'(z)\|_2 \leq O(1)\cdot\max\{\|X_i\|^2, \|Y_i\|^2\}\ .$$

The strong convexity of $f'(z)$ is given by

$$\nabla^2 f'(z) = B^{1/2}(\lambda I - M_s)B^{1/2} \succeq B^{1/2}(\frac{\delta_\times\sigma_k}{48}I)B^{1/2} \succeq \frac{\sigma_k\delta_\times\lambda_{\min}(B)}{48}\cdot I\ ,$$

where the first inequality is because (1) $\lambda_{\min}(\lambda I - M_s) \geq \frac{\delta_\times \lambda^*}{48}$ from Theorem 3.1, (2) $\lambda^* \geq \lambda_{\max}(M_s)$ from the definition of the algorithm, as well as (3) $\lambda_{\max}(M_s) \geq \sigma_k$ from the fact that $M_s$ is exactly $B^{-1/2}AB^{-1/2}$ but projecting out at most $k-1$ dimensions. $\quad\square$

**Claim F.5.** *With a preprocessing time $O(\mathsf{nnz}(X, Y))$, we can compute $(\nabla f_i'(z))w$ for an arbitrary vector $w$ and an arbitrary index $i \in [n]$ in time $O(kd)$, and $(\nabla f'(z))w$ in time $O(\mathsf{nnz}(X, Y) + kd)$.*

*Proof.* Note that $(\nabla f_i(z))w$ only costs running time $O(d)$ and thus it suffices to focus on the computation of $Qw$. We first note that one can pre-compute $A\mathbb{V}$, $B\mathbb{V}$ in time $O(\mathsf{nnz}(X, Y))$ with the help from the previous iteration. This is so because, suppose $\mathbb{V} = [\mathbb{V}', v]$ where $v$ is the new column vector added in the current iteration. Then, $A\mathbb{V} = [A\mathbb{V}', Av]$ and this additional column $Av$ can be computed in time $O(\mathsf{nnz}(X, Y))$. Next, we write

$$Q = (B\mathbb{V})(\mathbb{V}^\top A) + (A\mathbb{V})(\mathbb{V}^\top B) - (B\mathbb{V})(\mathbb{V}^\top A)\mathbb{V}(\mathbb{V}^\top B)$$

and every matrix between parenthesis of the above formulation has either only $O(k)$ rows or $O(k)$ columns. Therefore, computing $Qw$ costs a total running time at most $O(kd)$ and so is $(\nabla f_i'(z))w$.

Finally, we have $(\nabla f'(z))w = (\lambda B - A)w + Qw$ and therefore computing $(\nabla f'(z))w$ costs time $O(\mathsf{nnz}(X, Y) + kd)$. $\quad\square$

Since computing $(\nabla f_i'(z))w$ requires time $O(kd)$ and computing $(\nabla f'(z))w$ requires time $O(\mathsf{nnz}(X, Y) + kd)$, we can apply SVRG to minimize this sum-of-non-convex function $f'(z)$ [7, 24], with a running time

$$\widetilde{O}\Big(\big(\mathsf{nnz}(X, Y) + kd\big) + \frac{(\max\{\|X_i\|^2, \|Y_i\|^2\})^2}{\delta_\times^2 \sigma_k^2 \lambda_{\min}(B)^2} kd\big)\Big) = \widetilde{O}\Big(\mathsf{nnz}(X, Y) + \frac{(\kappa')^2}{\delta_\times^2 \sigma_k^2} kd\Big) \ .$$

Now, as long as $\mathsf{nnz}(X, Y) \geq kd$, we can apply acceleration on top of SVRG [7, 24] to given an accelerated running time

$$\widetilde{O}\Big(\mathsf{nnz}(X, Y) \cdot \Big(1 + \frac{\sqrt{\kappa'}}{\delta_\times^{1/2} \sigma_k^{1/2}(\mathsf{nnz}(X, Y)/(kd))^{1/4}}\Big)\Big) \ .$$

Taking into account all the iterations $s = 1, \ldots, k$ concludes the proof of Theorem F.3.

### F.4 Proofs of Theorems 6.2 and 6.3

The two corollaries follow from Theorem F.1 for the similar reason as Theorem 4.3 and Theorem 4.4 following from Theorem 4.1.

*Proof Theorem 6.2.* The approximation guarantees $\|\mathbb{V}_\phi^\top S_{xx} \mathbb{W}_\phi\|_2 \leq \varepsilon$ and $\|\mathbb{V}_\psi^\top S_{yy} \mathbb{W}_\psi\|_2 \leq \varepsilon$ follow from Theorem F.1.a, and the running time follows from Theorem F.2 and Theorem F.3. $\quad\square$

*Proof of Theorem 6.3.* The approximation guarantees $\max_{\phi \in \mathbb{R}^{d_x}, \psi \in \mathbb{R}^{d_y}} \big\{ \phi^\top S_{xy} \psi \ \big| \ \phi^\top S_{xx} \mathbb{V}_\phi = 0 \ \wedge$ $\psi^\top S_{yy} \mathbb{V}_\psi = 0 \big\} \leq (1+\varepsilon)\sigma_{k+1}$ follow from Theorem F.1.b and the definition of $M$, and the approximation guarantee $(1-\varepsilon)\sigma_i \leq |\phi_i' S_{xy} \psi_i| \leq (1+\varepsilon)\sigma_i$ follows from Theorem F.1.c and the fact that $|\lambda_{2i-1}| = |\lambda_{2i}| = \sigma_i$. The running time follows from Theorem F.2 and Theorem F.3. $\quad\square$